\documentclass[10pt]{amsart}

\usepackage{amsmath}
\usepackage{amsthm}
\usepackage{amssymb}

\input{diagrams}

\theoremstyle{plain}
\newtheorem{dfn}[subsection]{Definition}
\newtheorem{thm}[subsection]{Theorem}
\newtheorem{prp}[subsection]{Proposition}
\newtheorem{cor}[subsection]{Corollary}
\newtheorem{lma}[subsection]{Lemma}
\newtheorem{sdfn}[subsubsection]{Definition}

\newtheorem{sprp}[subsubsection]{Proposition}
\newtheorem{slma}[subsubsection]{Lemma}

\theoremstyle{remark}
\newtheorem{rmk}[subsection]{Remark}
\newtheorem{srmk}[subsubsection]{Remark}

\def\Spt{\mathbf{Spt}}
\def\Ho{\mathbf{Ho}}
\def\sg{\sigma}

\def\AA{\mathcal{A}}
\def\BB{\mathcal{B}}
\def\VV{\mathcal{V}}
\def\WW{\mathcal{W}}
\def\CC{\mathcal{C}}
\def\DD{\mathcal{D}}
\def\EE{\mathcal{E}}
\def\FF{\mathcal{F}}
\def\MM{\mathcal{M}}

\def\Sing{\mathrm{Sing}}
\def\Dop{\Delta^{\op}}
\def\lrto{\leftrightarrows}
\def\rlto{\rightleftarrows}
\def\ito{\rightarrowtail}
\def\eqv{\overset{\sim}{\longrightarrow}}
\def\cto{\twoheadrightarrow}
\def\Sets{\mathbf{Sets}}
\def\Cat{\mathbf{Cat}}
\def\GG{\mathbb{G}}
\def\RR{\mathbb{R}}
\def\inc{\hookrightarrow}
\def\Top{\mathbf{Top}}
\def\kk{\underline{k}}
\def\mm{\underline{m}}
\def\nn{\underline{n}}
\def\Sph{\mathbb{S}}
\def\Grph{\mathbf{Grph}}
\def\nGrph{\mathbf{nGrph}}

\def\nCat{\mathbf{nCat}}
\def\nplusCat{\mathbf{n\!+\!1Cat}}
\def\nplusGrph{\mathbf{n\!+\!1Grph}}
\def\NN{\mathcal{N}}
\def\BB{\mathcal{B}}
\def\dto{\rightrightarrows}

\def\lra{\longrightarrow}
\def\op{\mathrm{op}}
\def\CW{\mathrm{CW}}
\def\ZZ{\mathbb{Z}}

\begin{document}

\title[Iterated wreath-product of the simplex category]{Iterated wreath product of the simplex category and iterated loop spaces}

\author{Clemens Berger}

\date{9 December 2005; revised 21 June 2006.}
\begin{abstract}Generalising Segal's approach to $1$-fold loop spaces, the homotopy theory of $n$-fold loop spaces is shown to be equivalent to the homotopy theory of reduced $\Theta_n$-spaces, where $\Theta_n$ is an iterated wreath product of the simplex category $\Delta$. A sequence of functors from $\Theta_n$ to $\Gamma$ allows for an alternative description of the Segal spectrum associated to a $\Gamma$-space. 

In particular, each Eilenberg-MacLane space has a canonical reduced $\Theta_n$-set model. The number of $(n+d)$-dimensional cells of the resulting $\CW$-complex of type $K(\ZZ/2\ZZ,n)$ is a generalised Fibonacci number.\end{abstract}

\maketitle

\section*{Introduction}According to Segal \cite{Se}, reduced $\Gamma$-spaces are models for infinite loop spaces. His construction is based on an iterative use of reduced simplicial spaces as models for $1$-fold loop spaces. The purpose of the present text is to interpolate between these two constructions, and to show that the homotopy theory of $n$-fold loop spaces may be described inside the category of reduced $\Theta_n$-spaces, where $\Theta_n$ is a certain iterated wreath product of the simplex category $\Delta$. The \emph{duals} of these operator categories have been introduced by Joyal \cite{Joy1} as a possible starting point for the definition of weak higher categories. Therefore, our models for iterated loop spaces give some evidence to Baez and Dolan's \cite{BD} hypothesis that $n$-fold loop spaces are weak $n$-groupoids with one $0$-cell, one $1$-cell, one $2$-cell, ..., and one $(n-1)$-cell.
 
We begin our exposition with a possible characterisation of what it means to be a Quillen model category \cite{Q} for $n$-fold loop spaces, where $1\leq n\leq\infty$. In the extremal cases $n=1$ and $n=\infty$, this may be compared with Thomason and May's uniqueness results (cf. \cite{T} and \cite{MT}) for $1$-fold and infinite delooping machines. Our approach is based on a general definition of the \emph{derived image} of a right Quillen functor between model categories. The homotopy category of $n$-fold loop spaces (resp. infinite loop spaces) is then by definition the derived image of the $n$-fold loop functor $\Omega^n:\Top_*\to\Top_*$ (resp. of the right adjoint $\Omega^\infty:\Spt\to\Top_*$ of the suspension spectrum functor). As motivating examples, we show that the categories of reduced $\Delta$- and $\Gamma$-spaces nicely fit into our framework, insofar as they are model categories for the derived image of $\Omega^1$ and of $\Omega^\infty$ respectively.

Our main result (Theorem \ref{main}) states that the category of reduced $\Theta_n$-spaces is a model category for the derived image of $\Omega^n$. This statement is, roughly speaking, equivalent to the following property of reduced $\Theta_n$-spaces $X$ (cf. Lemma \ref{nfold}) : whenever $X$ is a cofibrant-fibrant reduced $\Theta_n$-space, the geometric realisation of $X$ is an \emph{$n$-fold delooping} of the underlying space of $X$.

Joyal \cite{Joy1} defines $\Theta_n$ to be the \emph{dual} of the category $\DD_n$ of finite combinatorial $n$-disks. Soon after Joyal's definition of $\DD_n$, Batanin and Street \cite{BS} have conjectured that $\Theta_n$ fully embeds in the category of strict $n$-categories by means of certain $n$-categorical analogs of the finite ordinals; this has been proved independently by Makkai and Zawadowsky \cite[Theorem 5.10]{MZ}, and the author \cite[Proposition 2.2]{mi}. The resulting description of the operators in $\Theta_n$ is however just as involved as Joyal's description of the dual operators in $\DD_n$, since it relies on Batanin's formula for the free $n$-category generated by certain $n$-graphs \cite{Ba1}, \cite[Definition 1.8]{mi}. We give here a new, conceptually simple, description of $\Theta_n$, together with a comparatively short proof of the duality between $\Theta_n$ and $\DD_n$. Indeed, $\Theta_n$ will be identified with an \emph{iterated wreath product} $\Delta\wr\cdots\wr\Delta$ of the simplex category $\Delta$. Many important properties of $\Theta_n$ follow now from the analogous properties of $\Delta$ by induction on $n$.

There is a formally similar wreath product for Segal's category $\Gamma$; moreover, disjoint sum induces a canonical functor $\Gamma\wr\Gamma\to\Gamma$. Therefore, Segal's functor $\gamma:\Delta\to\Gamma$ induces a whole sequence of \emph{assembly functors} $\gamma_n:\Theta_n\to\Gamma$, inductively defined as composite functors $\gamma_n:\Theta_n=\Delta\wr\Theta_{n-1}\overset{\gamma\wr\gamma_{n-1}}{\lra}\Gamma\wr\Gamma\to\Gamma.$ This assembly functor $\gamma_n$ also occurs (in dual form) in Batanin's $n$-operadic approach to $n$-fold loop spaces, cf. \cite{Ba2}, where it gives rise to an adjunction between \emph{$n$-operads} and \emph{symmetric operads}, the left adjoint of which sends contractible cofibrant $n$-operads to cofibrant $E_n$-operads, cf. Remark \ref{end}.

Most importantly for us, the assembly functors $\gamma_n$ relate the homotopy theories of reduced $\Gamma$-spaces and of reduced $\Theta_n$-spaces as follows : each reduced $\Gamma$-space $A$ determines an endofunctor $\underline{A}$ of the category of based spaces, together with a natural transformation of functors $\underline{A}(-)\wedge S^1\to\underline{A}(-\wedge S^1)$. The Segal spectrum of $A$ is by definition the spectrum $(\underline{A}(S^n))_{n\geq 0}$, where $S^n$ denotes the $n$-dimensional sphere. Whenever $A$ is a cofibrant-fibrant $\Gamma$-space with respect to Bousfield and Friedlander's \emph{stable} model structure \cite{BF}, the Segal spectrum of $A$ is an $\Omega$-spectrum, i.e. $\underline{A}(S^n)$ is an $n$-fold delooping of $\underline{A}(S^0)$. It turns out (see Corollary \ref{homeo}) that there is a functorial homeomorphism between $\underline{A}(S^n)$ and the geometric realisation of the inverse image $\gamma_n^*(A)$ of $A$; under these homeomorphisms, the structural maps of the Segal spectrum are induced by \emph{suspension} functors $\sg_n:\Theta_n\to\Theta_{n+1}$, already considered in \cite{Joy1}. The above-mentioned main property of reduced $\Theta_n$-spaces thus recovers Segal's construction of $\Omega$-spectra out of $\Gamma$-spaces.

Reduced $\Theta_n$-spaces may be compared with reduced $n$-simplicial spaces by means of an explicit functor $\Delta\times\cdots\times\Delta\to\Theta_n$. This diagonal functor induces a Quillen equivalence with respect to suitably chosen model structures (see Proposition \ref{compare1}) so that both reduced presheaf categories are model categories for $n$-fold loop spaces. The reduced $n$-simplicial model is obtained by a straightforward iteration of the Segal model for $1$-fold loop spaces, and it is folklore that this works out well, cf. Dunn \cite{Du} and Fiedorowicz-Vogt \cite{FV}. We construct also a Quillen equivalence (see Proposition \ref{compare2}) between reduced $\Theta_n$-spaces and $(n-1)$-reduced simplicial spaces. That the latter form a model category for $n$-fold loop spaces is the content of a widely circulated preprint of Bousfield \cite{Bou}.

There is a combinatorially interesting $\Theta_n$-set for each \emph{Eilenberg-MacLane space}, obtained from the $\Gamma$-set of the corresponding Eilenberg-MacLane spectrum by taking the inverse image with respect to $\gamma_n$. A counting lemma of Dolan \cite{Do} shows that the number of $(n+d)$-dimensional cells of the resulting $\CW$-complex of type $K(\ZZ/2\ZZ,n)$ is a generalised \emph{Fibonacci number}, the classical one's arising for $n=2$. For each finite abelian group $\pi$, the \emph{virtual Euler-Poincar\'e characteristic} of this reduced $\Theta_n$-set of type $K(\pi,n)$ is the expected one: the order of $\pi$ or its inverse depending on whether $n$ is even or odd.\vspace{1ex} 


The plan of this article is as follows:\vspace{1ex}

Section $1$ introduces the concept of the derived image of a right Quillen functor. We show that the derived image is uniquely determined up to equivalence of categories; moreover, we establish a rigidity result to the effect that any suitable Quillen adjunction between two models of the same derived image is a Quillen equivalence.\vspace{1ex}

Section $2$ discusses the Segal models for $1$-fold and infinite loop spaces from this model-theoretical viewpoint. The material of this section is classical, and we do not claim any originality for it. It is however pleasant to observe how nicely the existing literature fits into the language of derived image-factorisations.

Reduced simplicial spaces recently occured at several places (cf. Rezk \cite{R} and Bergner \cite{Ber}) in the context of the so-called \emph{Segal categories}. Our model structure on reduced simplicial spaces is a \emph{localisation} of the canonical injective model structure for reduced Segal categories; this reflects the fact that the fibrant objects of any model for the derived image of the loop functor have to be ``group-complete''.

Concerning $\Gamma$-spaces we follow as closely as possible the original texts of Segal \cite{Se} and Bousfield-Friedlander \cite{BF}; in particular, we use throughout the stable injective model structure on $\Gamma$-spaces as defined by Bousfield and Friedlander. This implies that our cofibrancy condition for $\Gamma$-spaces is somehow the weakest possible, but our fibrancy condition is stronger than Segal's condition of ``speciality''.\vspace{1ex}

Section $3$ is central : we introduce the wreath products over $\Delta$ and $\Gamma$, and study the $n$-fold wreath product $\Delta\wr\cdots\wr\Delta$, which serves as definition for $\Theta_n$. For consistency with the existing literature, we show that the so defined $\Theta_n$ densely embeds in the category of strict $n$-categories; Batanin's \cite{Ba1} star-construction for level-trees makes this embedding explicit. We then show that the presheaf topos on $\Theta_n$ is a classifying topos for Joyal's \cite{Joy1} combinatorial $n$-disks, and deduce the above-mentioned duality between $\Theta_n$ and $\DD_n$ from the fact that all idempotents in $\Theta_n$ split. In particular, this provides the category of $\Theta_n$-spaces with a left exact geometric realisation functor, already studied in \cite{Joy1} and \cite{mi}. We finally give an internal characterisation of the assembly functor $\gamma_n:\Theta_n\to\Gamma$ showing that the induced geometric morphism classifies the ``generic $n$-sphere'' in $\Theta_n$-sets.\vspace{1ex}

Section $4$ establishes the main result of this article, namely that the category of reduced $\Theta_n$-spaces is a model for the derived image of the $n$-fold loop functor. We compare reduced $\Theta_n$-spaces with reduced $n$-simplicial spaces and with $(n-1)$-reduced simplicial spaces, and show (using the rigidity result of Section $1$) that all three reduced presheaf categories carry Quillen equivalent model structures. We finally describe the reduced $\Theta_n$-set model for an Eilenberg-MacLane space of type $K(\pi,n)$ and determine its virtual Euler-Poincar\'e characteristic, as well as the relationship with generalised Fibonacci numbers.\vspace{1ex}

\emph{Acknowledgements.} The present text owes a lot to discussions I had with several friends and colleagues; I would like to thank especially Michael Batanin, Boris Chorny, Denis-Charles Cisinski, James Dolan, Bj\o rn Dundas, Hans-Werner Henn, Andr\'e Joyal, Jean-Louis Loday, Ieke Moerdijk, Ross Street and Rainer Vogt.


\section{The derived image of a right Quillen functor}

This section aims to give a sufficient set of axioms for a category to be called a \emph{model for the category of $n$-fold loop spaces}. In the literature, axiomatic appoaches exist for $1$-fold and infinite loop spaces, cf. Thomason \cite{T}, Fiedorowicz \cite{F}, May-Thomason \cite{MT}; our approach is based on the formalism of Quillen model structures. It turns out that we are looking for some kind of ``homotopical image-factorisation'' of the $n$-fold loop functor. In order to make this precise, we begin by defining what we mean by an honest image-factorisation of a functor. We apologise if ever the adopted terminology is not well suited from a purely categorical viewpoint.

We call \emph{left-conservative} (resp. \emph{right-conservative}) any isomorphism-reflecting left (resp. right) adjoint functor. The prototypical example of a left- (resp. right-) conservative functor is a comonadic (resp. monadic) functor.

We call \emph{reflection} (resp. \emph{coreflection}) any functor admitting a fully faithful right (resp. left) adjoint. Equivalently, this means that the counit (resp. unit) of the adjunction is an isomorphism. The prototypical example of a (co)reflection is the adjoint of the inclusion of a (co)reflective subcategory.

\begin{dfn}\label{imfact}An \emph{image-factorisation} of a right adjoint functor $G$ is a factorisation of $G$ into a coreflection followed by a right-conservative functor.

Dually, a \emph{coimage-factorisation} of a left adjoint functor $F$ is a factorisation of $F$ into a reflection followed by a left-conservative functor.\end{dfn}

Let $G:\EE'\to\EE$ be a right adjoint functor and let\begin{diagram}[noPS,small,p=1mm]&\EE'&\pile{\rTo^\Psi\\\lTo_\Phi}&\MM&\pile{\rTo^U\\\lTo_L}&\EE\end{diagram} be an image-factorisation of $G$. The intermediate category $\MM$ will be called an image of $G$. We shall usually abbreviate this image-factorisation by the triple $(\Phi,\MM,U)$. For two image-factorisations $(\Phi_1,\MM_1,U_1)$ and $(\Phi_2,\MM_2,U_2)$ of $G$, we say that an adjunction $d^*:\MM_1\lrto\MM_2:d_*$ is \emph{well-adapted} whenever $\Phi_1\cong\Phi_2d^*$ and $U_1d_*\cong U_2$. A well-adapted adjunction between coimages is defined dually.

\begin{lma}\label{imuniq}Any two images of a right adjoint (resp. coimages of a left adjoint) functor are equivalent as categories.\end{lma}

\begin{proof}Consider the following diagram of functors\begin{diagram}[noPS,small,p=1mm]\EE'&\pile{\rTo^{\Psi'}\\\lTo_{\Phi'}}&\MM'\\\dTo^\Psi\uTo_\Phi&&\dTo_{U'}\\\MM&\rTo_{U}&\EE\end{diagram}where the outer square commutes, $\Psi,\Psi'$ are coreflections, $U,U'$ are right-conservative, and $\Phi,\Phi'$ are left adjoint to $\Psi,\Psi'$. We have a chain of natural isomorphisms:$$U\Psi\Phi'\Psi'\Phi=U'\Psi'\Phi'\Psi'\Phi\cong U'\Psi'\Phi=U\Psi\Phi\cong U.$$Observe that the natural isomorphism $U'\eta'_{\Psi'\Phi}:U'\Psi'\Phi\cong U'\Psi'\Phi'\Psi'\Phi$ is inverse to the isomorphism $U'\Psi'\epsilon'_\Phi:U'\Psi'\Phi'\Psi'\Phi\cong U'\Psi'\Phi$, which may be identified with $U\Psi\epsilon'_\Phi:U\Psi\Phi'\Psi'\Phi\cong U\Psi\Phi$. Therefore, the composite isomorphism above belongs to the image of $U$, and since $U$ is conservative, this implies the existence of an isomorphism $\Psi\Phi'\Psi'\Phi\cong Id_\MM$. Similarly, $\Psi'\Phi\Psi\Phi'\cong Id_{\MM'}$. Therefore, $\Psi'\Phi$ and $\Psi\Phi'$ define an equivalence of categories between $\MM$ and $\MM'$. The proof of the dual statement is dual.\end{proof}

For the sequel, it will be important that the above \emph{uniqueness} result can be complemented by the following \emph{rigidity} result.

\begin{lma}\label{adjoints}Any well-adapted adjunction between two (co)images of the same functor is an adjoint equivalence.\end{lma}

\begin{proof}With the notations of the preceding proof, let $d^*:\MM\lrto\MM':d_*$ be a well-adapted adjunction. Then it follows from $U=U'd_*$ that $d_*$ is right-conservative; moreover, for each objet $X$ of $\MM$, the unit $X\to\Psi\Phi(X)$ factors as a composite of units $X\to d_*d^*(X)\to d_*\Psi'\Phi'd^*(X)$. Therefore, the right adjoint $d_*$ is not only conservative, but also a coreflection, and hence part of an adjoint equivalence with quasi-inverse $d^*$. The proof of the dual statement is dual.\end{proof}

\begin{rmk}(Co)image-factorisations appear naturally at different places in the literature; we shall briefly recall two of them: \emph{Beck's (co)monadicity theorem} and the \emph{``canonical'' factorisation} of a \emph{geometric morphism}. Indeed, let $F:\EE\lrto\EE':G$ be an adjunction, with unit $\eta:Id_\EE\to GF$, and counit $\epsilon:FG\to Id_{\EE'}$, and denote the category of algebras for the monad $GF$ by $\EE^{GF}$. There is a canonical comparison functor $\Psi:\EE'\to\EE^{GF}:Y\mapsto(GY,G\epsilon_Y)$; in particular, $G=U\Psi$ where $U:\EE^{GF}\to\EE$ is the forgetful functor. By definition, $G$ is monadic precisely when $\Psi$ is an equivalence of categories. The first of Beck's conditions for monadicity (namely that $G$ ``creates'' coequalisers for reflexive pairs in $\EE'$, for which the $G$-image has a split coequaliser in $\EE$) is \emph{equivalent} to the existence of a fully faithful left adjoint $\Phi$ of $\Psi$, i.e. to the condition that $G=U\Psi$ is an image-factorisation in the sense of Definition \ref{imfact}, cf. the dual of \cite[Lemma 1.1a-b]{AT}. The second condition for monadicity (namely that $G$ is right-conservative) is then equivalent to $\Psi$ being right-conservative, but this merely expresses the fact that a coreflection is an equivalence if and only if it is right-conservative. The discussion for Beck's comonadicity theorem is dual.

  A \emph{geometric morphism} between toposes is an adjoint pair $\phi_*:\EE'\rlto\EE:\phi^*$ such that the (left adjoint) inverse image functor $\phi^*$ is left exact (i.e. preserves finite limits). By definition, such a geometric morphism is an \emph{embedding} (resp. \emph{surjection}) if $\phi^*$ is a reflection (resp. faithful), cf. MacLane-Moerdijk \cite{MM}. Observe that $\phi^*$ is faithful if and only if it is conservative. Indeed, any faithful functor reflects monos and epis, and any topos is ``balanced'' (iso=mono+epi), so any faithful inverse image functor $\phi^*$ is conservative; conversely, any conservative left exact functor reflects finite limits, and so is faithful. Therefore, a geometric morphism $\phi_*:\EE'\rlto\EE:\phi^*$ factors as surjection followed by an embedding if and only if $\phi^*$  factors as a reflection followed by a left-conservative functor, i.e. iff $\phi^*$ admits a coimage-factorisation in the sense of Definition \ref{imfact}. The existence of such a coimage-factorisation follows from Beck's comonadicity theorem like above, since the left exactness of $\phi^*$ and the finite completeness of toposes give the first of Beck's conditions for free. The essential uniqueness can be deduced from Lemma \ref{adjoints}.

General criteria for the existence of a coimage-factorisation have been given by Applegate-Tierney \cite{AT} and Day \cite{Day}. The dual situation of an image-factorisation has been studied by Adamek-Herrlich-Tholen \cite{AHT}.\end{rmk}

  We shall now define the homotopical analog of (co)image-factorisations in the framework of Quillen's closed model categories. For the convenience of the reader unfamiliar with this theory, we give a short r\'esum\'e of the main definitions. Excellent references on the subject include \cite{Q}, \cite{DHK}, \cite{Ho}, \cite{Hir}. As is nowadays usually the case, we shall omit the adjective `closed'. 

A \emph{model category} is a finitely complete and finitely cocomplete category, equipped with three distinguished classes of morphisms, called respectively \emph{cofibrations}, \emph{weak equivalences} and \emph{fibrations}. Usually, these morphisms are depicted by arrows of the form $\ito,\eqv,\cto$; the following four axioms have to be satisfied:\vspace{1ex}

\begin{itemize}\item[(M1)]If any two among $f,g$ and $gf$ are weak equivalences, then so is the third;\item[(M2)]Cofibrations, weak equivalences and fibrations compose, contain all isomorphisms and are closed under retract;\item[(M3)]Cofibrations (resp. trivial cofibrations) have the left lifting property with respect to trivial fibrations (resp. fibrations);\item[(M4)]Any morphism factors as a cofibration followed by a trivial fibration, and as a trivial cofibration followed by a fibration.\end{itemize}
Here, a trivial (co)fibration means a (co)fibration which is also a weak equivalence. An object $X$ of $\EE$ is called \emph{cofibrant} (resp. \emph{fibrant}) if the unique morphism from an inital object of $\EE$ to $X$ (resp. from $X$ to a terminal object of $\EE$) is a cofibration (resp. fibration). For a general object $X$, a \emph{cofibrant replacement} consists of a cofibrant object $c_\EE(X)$ together with a weak equivalence $c_\EE(X)\eqv X$, while a \emph{fibrant replacement} consists of a fibrant object $f_\EE(X)$ together with a weak equivalence $X\eqv f_\EE(X)$. Axiom (M4) implies the existence of cofibrant and fibrant replacements for any object $X$. 

A \emph{fibrant} (resp. \emph{cofibrant}) \emph{replacement functor} for $\EE$ is an endofunctor $f_\EE$ (resp. $c_\EE$) of $\EE$ endowed with a natural transformation $id_\EE\to f_\EE$ (resp. $c_\EE\to id_\EE$) which is pointwise a fibrant (resp. cofibrant) replacement. The existence of such replacement functors is not a formal consequence of Quillen's axioms; one way to obtain such replacement functors is to require functoriality of the factorisations in (M4); functorial factorisations exist for instance in any \emph{cofibrantly generated} model category, which will always be the case for us, cf. Section \ref{cofgen}.

There is a well defined \emph{homotopy relation} $\sim$ on the set $\EE(X,Y)$ of morphisms from $X$ to $Y$, whenever $X$ is cofibrant and $Y$ is fibrant. This leads to the following definition of the \emph{homotopy category} $\Ho(\EE)$: both categories $\EE$ and $\Ho(\EE)$ have the same objects, and for each pair of objects, we have  $$\Ho(\EE)(X,Y)=\EE(c_\EE(X),f_\EE(Y))/\sim.$$ This definition does not depend (up to canonical bijection) on the choice of the replacements; moreover, homotopy classes compose in such a way that passage to the homotopy class defines a functor $\gamma_\EE:\EE\to\Ho(\EE)$. The functor $\gamma_\EE$ is \emph{initial among functors out of $\EE$ turning weak equivalences into isomorphisms}; moreover, a morphism in $\EE$ is a weak equivalence if and only if its image under $\gamma_\EE$ is an isomorphism in $\Ho(\EE)$;  The homotopy category $\Ho(\EE)$ is thus an explicit model of the ``localisation'' of $\EE$ with respect to the class of weak equivalences, cf. Gabriel and Zisman \cite{GZ}.

A \emph{Quillen adjunction} $F:\EE\lrto{\EE'}:G$ between model categories is an adjunction between the underlying categories with the property that the left adjoint $F$ preserves cofibrations and the right adjoint $G$ preserves fibrations. This implies (cf. \cite[7.7/8.5]{Hir}) that $F$ also preserves trivial cofibrations as well as weak equivalences between cofibrant objects, and $G$ also preserves trivial fibrations as well as weak equivalences between fibrant objects. Usually, the left (resp. right) adjoint of a Quillen adjunction is called a \emph{left} (resp. \emph{right}) \emph{Quillen functor}.

The \emph{left derived functor} $LF:\Ho(\EE)\to\Ho({\EE'})$ is induced by the universal property of $\Ho(\EE)$ applied to the composite functor $\gamma_{\EE'}\circ F\circ c_\EE$, where $c_\EE$ is a cofibrant replacement functor. The \emph{right derived functor} $RG:\Ho({\EE'})\to\Ho(\EE)$ is induced by the universal property of $\Ho(\EE')$ applied to $\gamma_\EE\circ G\circ f_{\EE'}$, where $f_{\EE'}$ is a fibrant replacement functor. The derived functors of a Quillen adjunction define a \emph{derived adjunction} $LF:\Ho(\EE)\lrto\Ho({\EE'}):RG$ between the homotopy categories. A Quillen adjunction is called a \emph{Quillen equivalence} if this derived adjunction is an equivalence of categories.

The derived adjunction may be represented by means of the following binatural bijections, where $X$ is a cofibrant object of $\EE$, and $Y$ a fibrant object of $\EE'$:$$\Ho(\EE)((LF)(X),Y)\cong\EE(F(X),Y)/\sim\,\cong\EE(X,G(Y))/\sim\,\cong\Ho(\EE)(X,(RG)(Y)).$$The unit at $X$ of the derived adjunction is thus represented by the adjoint of any fibrant replacement $F(X)\eqv f_{\EE'}(F(X))$; this adjoint $X\to G(f_{\EE'}(F(X)))$ will be called a \emph{homotopy-unit}; dually, the adjoint $F(c_\EE(G(Y)))\to Y$ of any cofibrant replacement $c_\EE(G(Y))\eqv G(Y)$ will be called a \emph{homotopy-counit}, since it represents the counit at $Y$ of the derived adjunction. Observe that homotopy (co)units are only defined for cofibrant objects of $\EE$ and fibrant objects of $\EE'$; nervertheless, the homotopy-(co)units represent up to isomorphism all (co)units of the derived adjunction, since the cofibrant objects of $\EE$ and fibrant objects of $\EE'$ represent up to isomorphism all objects of $\Ho(\EE)$ and $\Ho(\EE')$.

A left Quillen functor $F:\EE\to\EE'$ will be called \emph{homotopy-left-conservative} whenever a morphism $f$ between cofibrant objects of $\EE$ is a weak equivalence in $\EE$ if and only if $F(f)$ is a weak equivalence in $\EE'$. Dually, a right Quillen functor $G:\EE'\to\EE$ will be called \emph{homotopy-right-conservative} whenever a morphism $g$ between fibrant objects of $\EE'$ is a weak equivalence in $\EE'$ if and only if $G(g)$ is a weak equivalence in $\EE$. 

A left Quillen functor will be called a \emph{homotopy-reflection} if all homotopy-counits of the Quillen adjunction are weak equivalences. Dually, a right Quillen functor will be called a \emph{homotopy-coreflection} if all homotopy-units of the Quillen adjunction are weak equivalences. Dugger \cite{Du} uses the more suggestive term \emph{homotopy-surjection} for homotopy-reflections. The preceding considerations immediately imply

\begin{lma}\label{hom}A left Quillen functor is homotopy-left-conservative (resp. a homoto-py-reflection) if and only if its left derived functor is left-conservative (resp. a reflection). A right Quillen functor is homotopy-right-conservative (resp. a homotopy-coreflection) if and only if its right derived functor is right-conservative (resp. a coreflection).\end{lma}

The following criteria for a Quillen adjunction to be a Quillen equivalence are just homotopy versions of well known criteria for an adjunction to be an equivalence.

\begin{lma}\label{Qequiv}For a Quillen adjunction $F:\EE\lrto\EE':G$ the following five conditions are equivalent, cf. \cite[1.3.12/13/16]{Ho} :\begin{itemize}\item[(a)]The Quillen adjunction is a Quillen equivalence;\item[(b)]For any cofibrant object $X$ of $\EE$ and any fibrant object $Y$ of $\EE'$, a morphism $F(X)\to Y$ is a weak equivalence in $\EE'$ if and only if the adjoint morphism $X\to G(Y)$ is a weak equivalence in $\EE$;\item[(c)]All homotopy-units and homotopy-counits are weak equivalences;\item[(d)]$F$ is a homotopy-left-conservative homotopy-reflection;\item[(e)]$G$ is a homotopy-right-conservative homotopy-coreflection.\end{itemize}\end{lma}

\begin{dfn}A \emph{derived image} of a right Quillen functor $G:\EE'\to\EE$ is the homotopy category $\Ho(\MM)$ associated to any factorisation of $G$ into a homotopy-coreflection $\EE'\to\MM$ followed by a homotopy-right-conservative functor $\MM\to\EE$. A \emph{derived coimage} of a left Quillen functor is defined dually.\end{dfn}

Lemmas \ref{imuniq} and \ref{hom} immediately imply the following \emph{uniqueness} result:

\begin{prp}\label{uniq}Any two derived images (resp. coimages) of a right (resp. left) Quillen functor are equivalent as categories.\end{prp}

This uniqueness result is not optimal insofar as there is no indication whether the stated equivalence of homotopy categories is induced by a Quillen equivalence, or at least by a finite chain of Quillen equivalences. It is however unclear whether such a stronger uniqueness holds in general without further constraint on the factorisation.\vspace{1ex}
 
A model category $\MM$ like in the definition above will also be called a \emph{model} for the \emph{derived image} of $G$. Such a model is thus given by a triple $(\Phi,\MM,U)$ such that $\Phi:\MM\to\EE$ has a right adjoint homotopy-coreflection $\Psi:\EE\to\MM$, the functor $U:\MM\to\EE'$ is homotopy-right-conservative, and the composite functor $U\Psi$ is $G$.

A Quillen adjunction $d^*:\MM_1\lrto\MM_2:d_*$ between two models $(\Phi_1,\MM_1,U_1)$ and $(\Phi_2,\MM_2,U_2)$ of the derived image of $G$ will be called \emph{well-adapted} if $\Phi_1\cong\Phi_2d^*$ and $U_1d_*\cong U_2$.\vspace{1ex}

Lemmas \ref{adjoints} and \ref{hom} immediately imply the following \emph{rigidity} result:

\begin{prp}\label{rigid}Any well-adapted Quillen adjunction between models for the derived image of the same functor is a Quillen equivalence.\end{prp}

We shall say that $(\Phi,\MM,U)$ is a \emph{good model} whenever $\Phi$ \emph{preserves cofibrant-fibrant objects}; this is for instance always the case if every object of $\EE$ is fibrant. Good models have the advantage over general models that the homotopy-unit condition can be rephrased in a nice way; namely, for a good model $(\Phi,\MM,U)$, a homotopy-unit at a cofibrant object $X$ of $\MM$ is given by the composite $$X\to f_\MM(X)\to \Psi\Phi(f_\MM(X))$$ of a fibrant replacement in $\MM$ with the ordinary unit of the adjunction. Therefore, the condition that the homotopy-units are weak equivalences is equivalent (for good models) to the condition that the ordinary units are weak equivalences at cofibrant-fibrant objects. Since $U$ is homotopy-right-conservative, the latter condition is equivalent to the condition that the canonical map $$U(X)\to G\Phi(X),$$induced by the unit of the $\Phi$-$\Psi$-adjunction, is a weak equivalence for each cofibrant-fibrant object $X$ of $\MM$. In practice, it is this last condition which is most easily verified.

Explicitly, a model for the derived image of the $n$-fold loop functor (for short: a \emph{model for $n$-fold loop spaces}) consists of a triangle of Quillen adjunctions\begin{diagram}[noPS,p=1mm]\Top_*&\pile{\lTo^\Phi\\\rTo_\Psi}&\MM\\&\pile{\rdTo^{\Omega^n}\\\luTo_{\Sigma^n}}&\dTo^U\uTo_L\\&&\Top_*\end{diagram}with commuting inner (resp. outer) triangle of right (resp. left adjoints) such that $U$ is homotopy-right-conservative, and such that the homotopy-units of the $\Phi$-$\Psi$-adjunction are weak equivalences. We shall call $U$ the \emph{underlying-space functor} and $\Phi$ the \emph{$n$-fold delooping functor} of $\MM$.

On the category $\Top$ of compactly generated spaces, we use Quillen's model structure \cite{Q} with \emph{weak homotopy equivalences} as weak equivalences and \emph{Serre fibrations} as fibrations. The class of cofibrations contains the \emph{relative $CW$-complexes}. An explicit generating set for cofibrations (resp. trivial cofibrations) is the set of sphere-inclusions $S^{n-1}\inc B^n,\, n\geq 0,$ (resp. of ball-inclusions $B^{n-1}\inc B^n,\, n> 0$). The category $\Top_*$ of based objects in $\Top$ inherits a cofibrantly generated model structure from $\Top$ for which the $\Sigma^n$-$\Omega^n$ adjunction is a Quillen adjunction. Indeed, this model structure is obtained from Proposition \ref{transfer} applied to the adjunction $(-)_+:\Top\lrto\Top_*:U$, the left adjoint of which simply adds a disjoint base-point. In particular, the generating (trivial) cofibrations for the model structure on $\Top_*$ are the sets $\{S_+^{n-1}\inc B^n_+|\,n\geq 0\}$ and $\{B_+^{n-1}\inc B^n_+|\,n>0\}$.

Since every object is fibrant in $\Top_*$, any model is good; therefore,\begin{lma}\label{nfold}With the notations above, $(\Phi,\MM,U)$ is a model for $n$-fold loop spaces if and only if the underlying-space functor $U:\MM\to\Top_*$ is homotopy-right-conservative, and  the unit of the $\Phi$-$\Psi$-adjunction induces a weak equivalence $U(X)\eqv\Omega^n\Phi(X)$ for each cofibrant-fibrant object $X$ of $\MM$.\end{lma}

Similarly, a model for the derived image of $\Omega^\infty$ (for short: a \emph{model for infinite loop spaces}) consists of a triangle of Quillen adjunctions\begin{diagram}[noPS,silent,p=1mm]\Spt&\pile{\lTo^\Phi\\\rTo_\Psi}&\MM\\&\pile{\rdTo^{\Omega^\infty}\\\luTo_{\Sigma^\infty}}&\dTo^U\uTo_L\\&&\Top_*\end{diagram}
with commuting inner (resp. outer) triangle of right (resp. left adjoints) such that $U$ is homotopy-right-conservative, and such that the homotopy-units of the $\Phi$-$\Psi$-adjunction are weak equivalences. We shall call $U$ the \emph{underlying-space functor} and $\Phi$ the \emph{delooping functor} of $\MM$.

We use the stable model category $\Spt$ of topological spectra as defined by Bousfield and Friedlander \cite[2.1/5]{BF}. In particular, a spectrum $E=(E_n)_{n\geq 0}$ means a sequence of pointed compactly generated spaces $E_n$, equipped with structural maps $E_n\wedge S^1\to E_{n+1}$. The suspension spectrum functor $\Sigma^\infty:\Top_*\to\Spt$ is defined by $(\Sigma^\infty X)_n=X\wedge S^n$ with the obvious structural maps; its right adjoint $\Omega^\infty$ is given by $\Omega^\infty E=E_0$. 

A spectrum $E$ is called an \emph{$\Omega$-spectrum} if the adjoints $E_n\to\Omega E_{n+1}$ of the structural maps are weak equivalences. There is an endofunctor $T:\Spt\to\Spt$ converting spectra into $\Omega$-spectra; explicitly, $T(E)_n=\varinjlim_k\Omega^kE_{n+k}$. A map of spectra $f:E\to F$ is called a \emph{stable equivalence} if the induced map of $\Omega$-spectra $T(f):T(E)\to T(F)$ yields a weak equivalence $T(f)_n:T(E)_n\to T(F)_n$ in $\Top_*$ for each $n\geq 0$. A map of spectra $f:E\to F$ is called a \emph{stable fibration} if, for each $n\geq 0$, $f_n:E_n\to F_n$ is a fibration in $\Top_*$, and the induced map $E_n\to \Omega E_{n+1}\times_{\Omega F_{n+1}}F_n$ is a weak equivalence in $\Top_*$. The stable model structure on $\Spt$ has the stable equivalences as weak equivalences and the stable fibrations as fibrations. Our description of the stable fibrations differs from the original one, but Schwede \cite[A3]{Sch} shows that the two definitions are equivalent, and that the stable model category comes equipped with canonical generating sets for cofibrations and trivial cofibrations. The description of the latter also implies that the $\Sigma^\infty$-$\Omega^\infty$-adjunction is a Quillen adjunction.

Since in the stable model category $\Spt$, the fibrant objects are precisely the $\Omega$-spectra, we have to be a little more careful in analysing homotopy-units than for $n$-fold loop spaces. However, we shall see that in the case of interest for us, the delooping functor $\Phi$ sends cofibrant-fibrant objects to $\Omega$-spectra. Thus, it is enough to state

\begin{lma}\label{infinite}With the notations above, $(\Phi,\MM,U)$ is a good model for infinite loop spaces if and only if the underlying-space functor $U:\MM\to\Top_*$ is homotopy-right-conservative, the delooping functor $\Phi:\MM\to\Spt$ sends cofibrant-fibrant objects to $\Omega$-spectra, and the unit of the $\Phi$-$\Psi$-adjunction induces a weak equivalence $U(X)\eqv\Omega^\infty\Phi(X)$ for each cofibrant-fibrant object $X$ of $\MM$.\end{lma}

\subsection{Transfer and localisation of cofibrantly generated model structures}\label{cofgen}
We briefly recall two methods of defining model structures out of existing ones: \emph{transfer} and \emph{localisation}; a thorough discussion may be found in the book of Hirschhorn, to which we also refer for the omitted proofs, cf. \cite[Thms. 11.3.1-2]{Hir}.

The important concept is that of a \emph{cofibrantly generated} model category $\EE$. This means that $\EE$ is cocomplete, and that there exists a \emph{set} $I$ (resp. $J$) of cofibrations (resp. trivial cofibrations) which determines the class of trivial fibrations (resp. fibrations) by the appropriate lifting condition. Moreover, these sets are assumed to permit Quillen's \emph{small object argument}, i.e., the domains of the morphisms in $I$ (resp. $J$) have to be \emph{small} with respect to \emph{sequential} colimits of morphisms obtained from $I$ (resp. $J$) by arbitrary pushout. 

\begin{sprp}\label{transfer}Let $F:\EE\lrto\EE':G$ be an adjunction with left adjoint $F$ and right adjoint $G$. Assume that $\EE$ is a cofibrantly generated model category with generating sets $I$ and $J$, and that $\EE'$ is cocomplete and finitely complete.

If the sets $F(I)$ and $F(J)$ permit Quillen's small object argument and $G$ takes any sequential colimit of pushouts of morphisms in $F(J)$ to a weak equivalence in $\EE$, then $\EE'$ has a model structure, cofibrantly generated by $F(I)$ and $F(J)$, for which the weak equivalences (resp. fibrations) are the morphisms sent to weak equivalences (resp. fibrations) in $\EE$ under $G$.\end{sprp}

\begin{srmk}\label{specialtransfer}The hypothesis for transfer is in  particular fulfilled if $G$ preserves colimits, and $GF$ preserves cofibrations and trivial cofibrations. \end{srmk}

\begin{sdfn}\label{complete}A \emph{localising system} $L$ for a model category $\EE$ consists of a class of morphisms $L$, called \emph{locally trivial} or \emph{local equivalences}, such that\begin{enumerate}\item any weak equivalence is a local equivalence;\item local equivalences satisfy axioms (M1) and (M2) of a model category;\item locally trivial cofibrations are closed under pushout and sequential colimit.\end{enumerate}
A localising system $L$ is \emph{complete} if there exists a \emph{set} $S$ of locally trivial cofibrations such that\begin{enumerate}\item[(4)]the domains of the elements of $S$ are small with respect to sequential colimits of cofibrations;\item[(5)] any locally trivial $S$-fibration is a trivial fibration, where $S$-fibrations are those fibrations which have the right lifting property with respect to $S$.\end{enumerate}\end{sdfn}

\begin{sprp}\label{localisation}Let $\EE$ be a cofibrantly generated model category with complete localising system $(L,S)$. There exists a cofibrantly generated model category $\EE_S$ with same underlying category and same class of cofibrations as $\EE$, the so-called \emph{left Bousfield localisation of $\EE$ with respect to $L$}, whose weak equivalences are the local equivalences, and whose fibrations are the $S$-fibrations. 

In particular, the class of local equivalences coincides with the class of so-called $S$-local equivalences, i.e. those $f:X\to Y$ which induce a bijection $\Ho_\EE(f,Z):\Ho_\EE(Y,Z)\eqv\Ho_\EE(X,Z)$ for each $S$-fibrant object $Z$.\end{sprp}

\begin{srmk}\label{local}Since $\EE$ and $\EE_S$ have the same class of cofibrations, the identity functor $\EE\to\EE_S$ is a left Quillen functor and the identity functor $\EE_S\to\EE$ is a right Quillen functor. Moreover, morphisms \emph{between $S$-fibrant objects} are local equivalences (resp. $S$-fibrations) if and only they are weak equivalences (resp. fibrations). In particular, the identity functor $\EE_S\to\EE$ is \emph{homotopy-right-conservative}.

The definition of an $S$-local equivalence only depends on $S$ and the initial model category $\EE$ so that the axioms of a complete localising sytem $(L,S)$ could have been phrased entirely in terms of $S$. Indeed, important existence theorems of Bousfield, Smith and Hirschhorn give general conditions under which the class $L_S$ of $S$-local equivalences for \emph{any} given set of cofibrations $S$ is part of a complete localising system. For instance, if $\EE$ is left proper and \emph{cellular} or \emph{accessible} (see Hirschhorn \cite{Hir} for the precise definitions), then there is a suitable ``saturation'' $\bar{S}$ of $S$, for which $(L_S,\bar{S})$ is a complete localising system in the sense of Definition \ref{complete}.

In the literature, Proposition \ref{localisation} is often used without explicit mention. For instance, the \emph{stable} model structure on $\Spt$ is obtained this way by localisation of a \emph{strict} model structure with levelwise weak equivalences, cf. Schwede \cite[A3]{Sch}. Other examples may be found in \cite{RSS}, or in \cite[chapter 4]{mi}.\end{srmk}

\section{The Segal models for $1$-fold and infinite loop spaces}

In this section we present the Segal models for $1$-fold and infinite loop spaces (namely reduced $\Delta$- and reduced $\Gamma$-spaces) using the language of Section $1$. 

\subsection{Reduced simplicial spaces}\label{redsimpspaces}
The \emph{simplex category} $\Delta$ is the category of finite non-empty ordinals $[m]=\{0,\dots,m\},\,m\geq 0,$ and order-preserving maps. There is a well known embedding $\Delta^-:\Delta\inc\Top$ sending the ordinal $[m]$ to the standard euclidean $m$-simplex $\Delta^m$. The induced realisation functor $\Top^{\Dop}\to\Top$, defined by left Kan extension $X\mapsto X\otimes_\Delta\Delta^-$, will be denoted by the symbol $|\!-\!|_\Delta$. 

A \emph{simplicial space} $X:\Dop\to\Top$ is \emph{reduced} if $X([0])$ is the one-point space. The category of reduced simplicial spaces is denoted by $\Top_{red}^{\Dop}$. The right adjoint of $|\!-\!|_{\Delta}:\Top_{red}^{\Dop}\to\Top_*$ is the \emph{Segal loop functor} $\Omega_{Seg}:\Top_*\to\Top_{red}^{\Dop}$, defined by $\Omega_{Seg}(X)([m])=\Top((\Delta_m,sk_0\Delta_m),(X,*))$ with the usual topology; in particular, $\Omega_{Seg}(X)([1])=\Omega(X)$. The \emph{underlying space} of a reduced simplicial space is $U(Y)=Y([1])$. This defines the following commutative triangle of adjunctions:\begin{gather}\label{1folddiag}\begin{diagram}[noPS,p=1mm]\Top_*&\pile{\lTo^{|\!-\!|_{\Delta}}\\\rTo_{\Omega_{Seg}}}&\Top_{red}^{\Dop}\\&\pile{\rdTo^{\Omega}\\\luTo_{\Sigma}}&\dTo^U\uTo_L\\&&\Top_*\end{diagram}\end{gather}In order to specify a class of reduced simplicial spaces for which the realisation functor is a delooping, Segal \cite{Se} introduces the following three conditions:\begin{itemize}\item[(a)]For each $m\geq 1$, the canonical map $X([m])\to X([1])\times_{X([0])}\cdots\times_{X([0])} X([1])$ induced by the ``outer face'' operators $[0]\dto[1]\to[m]$ is a weak equivalence;\item[(b)]$X$ is ``group-like'';\item[(c)]$X$ is ``good''.\end{itemize}A simplicial operator $\phi:[k]\to[m]$ is an \emph{outer} face operator if $\phi(i+1)=\phi(i)+1$ for $i=0,\dots,k-1$. Condition (a) implies that the set of path-components $\pi_0(X([1]))$ carries a canonical monoid structure. By condition (b), this monoid has to be a group. Condition (c) means that $X$ belongs to a class of simplicial spaces, for which the realisation functor $|\!-\!|_\Delta$ takes pointwise weak equivalences to weak equivalences. 

Segal's key proposition \cite[Proposition 1.5]{Se} states that for reduced simplicial spaces $X$ satisfying (a), (b), (c), the canonical map $U(X)\to\Omega(|X|_\Delta)$ is a weak equivalence. Therefore, according to Lemma \ref{nfold}, in order to show that reduced simplicial spaces form a model category for $1$-fold loop spaces, it suffices to find a model structure for which $U$ is homotopy-right-conservative, and for which the cofibrant-fibrant objects satisfy conditions (a), (b), (c) above. Theorem \ref{Segalfibrant} below provides such a model structure and actually reproves to some extent Segal's key proposition.

Since $\Delta$ is a \emph{Reedy category}, cf. \cite{Hir}, \cite{Ho} and Section \ref{factor}, there is a canonical Reedy model structure on simplicial spaces: the weak equivalences are the pointwise weak equivalences, and the fibrations are the Reedy-fibrations, i.e. those natural transformations $X\to Y$ for which the so-called ``matching'' maps $$X([m])\to Y([m])\times_{Y(\partial\Delta[m])}X(\partial\Delta[m])$$ are fibrations in $\Top$ for all $[m]$. As usual, $\partial\Delta[m]$ denotes the boundary of the representable presheaf $\Delta[m]=\Delta(-,[m])$ (i.e. the union of all proper faces of $\Delta[m]$), and $X(\partial\Delta[m])$ is shorthand notation for the inverse limit $\varprojlim_{([k]\to[m])\in\partial\Delta[m]}X([k])$.

A morphism of simplicial spaces $f:X\to Y$ will be called a \emph{realisation weak equivalence} if the left derived functor of the realisation functor takes $f$ to an isomorphism in $\Ho(\Top)$. 

Following \cite[Definition 3.3]{RSS}, a Reedy-fibration of simplicial spaces $f:X\to Y$ will be called \emph{equifibered} if, for each injective simplicial operator $[k]\to[m]$, the induced map $X([m])\to Y([m])\times_{Y([k])}X([k])$ is a weak equivalence. 

The following theorem is a special case of Theorem 3.6(1) of Rezk, Schwede and Shipley \cite{RSS}. It is a quite direct application of Proposition \ref{localisation} and uses in an essential way a fundamental lemma of V. Puppe concerning the compatibility of the realisation functor with homotopy pullbacks, cf. \cite[Proposition 4.4]{RSS}.

\begin{thm}\label{equifibered}The Reedy model category of simplicial spaces admits a left Bousfield localisation with respect to realisation weak equivalences, the fibrations of which are precisely the equifibered Reedy-fibrations.\end{thm}

\begin{dfn}\label{Segalfibration}A Reedy-fibration of simplicial spaces $f:X\to Y$ is called a \emph{Segal-fibration} if the following two conditions hold:\begin{itemize}\item[(a)]for each $m\geq 1$, the induced map $X([m])\to Y([m])\times_{Y(G_m)}X(G_m)$ is a weak equivalence, where $G_m$ is the simplicial subset of $\Delta[m]$ defined by$$G_m=\bigcup_{\phi\in\Delta([1],[m])\,|\,\phi(1)=\phi(0)+1}\phi(\Delta[1]);$$
\item[(b)]the induced map $PX\to PY\times_YX$ is an \emph{equifibered} Reedy-fibration of simplicial spaces, where $PX\to X$ (resp. $PY\to Y$) denotes the canonical ``path-fibration'' of $X$ (resp. $Y$).\end{itemize}\end{dfn}
The canonical ``path-fibration'' is defined as follows: $(PX)([n])=X([n+1])$ with simplicial structure given by forgetting the last face- and degeneracy operators of $X$. These last face operators induce a morphism $PX\to X$ (which is however a Reedy-fibration only for Reedy-fibrant $X$), while the last degeneracy operators induce a combinatorial deformation-retraction of $PX$ onto $X([0])$.

For any reduced simplicial space $X$, conditions (a) and (b) for a Segal-fibration $X\to *$ translate into Segal's conditions (a) and (b). In particular, a reduced Reedy-fibrant simplicial space $X$, which satisfies Segal's condition (a), is group-like \emph{if and only if} its path-fibration $PX\to X$ is equifibered; this fact has already been exploited by Segal, cf. his proof of \cite[Proposition 1.5-6]{Se}.

For any Reedy-fibration $X\to Y$, the induced map $PX\to PY\times_YX$ is again a Reedy-fibration; if the former is equifibered, then so is the latter; therefore, condition (b) of a Segal-fibration is automatically fulfilled for equifibered Reedy-fibrations; however, condition (b) of a Segal-fibration $f$ is much weaker than $f$ being an equifibered Reedy-fibration. For instance, an equifibered Reedy-fibration between reduced simplicial spaces is necessarily a trivial Reedy-fibration, which is not the case for Segal-fibrations.

\begin{thm}\label{Segalfibrant}Pointwise weak equivalences and Reedy-fibrations induce a cofibrantly generated model structure on \emph{reduced} simplicial spaces. This model category admits a left Bousfield localisation with respect to realisation weak equivalences, the fibrations of which are precisely the \emph{Segal-fibrations}. 

The localised model category of reduced simplicial spaces is a model for $1$-fold loop spaces.\end{thm}

\begin{proof}The Reedy model structure on $\Top_*^{\Dop}$ induces by transfer a Reedy model structure on $\Top_{red}^{\Dop}$ since reduced simplicial spaces form a bireflective subcategory of based simplicial spaces, and Remark \ref{specialtransfer} applies. The class of realisation weak equivalences is a localising system in the sense of Definition \ref{complete}. Moreover, Segal-fibrations can be characterised (among Reedy-fibrations) by a right lifting property with respect to a set $S$ of realisation trivial cofibrations which satisfy condition (4) of \ref{complete} (this is clear for condition (a) of a Segal-fibration, and requires some combinatorics for condition (b)).

In order to apply Proposition \ref{localisation}, it remains to be shown that realisation-trivial Segal-fibrations $f:X\to Y$ are trivial Reedy-fibrations. The induced map $\Omega f:\Omega X\to\Omega Y$ on the fiber of the ``path-fibrations'' is an equifibered Reedy-fibration by condition (b) of a Segal-fibration. Since $f$ is a realisation weak equivalence, and the realisation functor is left exact, the induced map $PX\to PY\times_YX$ is a realisation-trivial equifibered Reedy-fibration, and hence a trivial Reedy-fibration by Theorem \ref{equifibered}. Therefore, $\Omega f$ is a trivial Reedy-fibration, and $X([1])=(\Omega X)([0])\to(\Omega Y)([0])=Y([1])$ is a trivial fibration in $\Top$. By condition (a) of a Segal-fibration, the commutative square\begin{gather}\label{Segal}\begin{diagram}[noPS,small]X([m])&\rTo&Y([m])\\\dTo&&\dTo\\X([1])\times\cdots\times X([1])&\rTo&Y([1])\times\cdots\times Y([1])\end{diagram}\end{gather}is homotopy cartesian, so that $f$ is a pointwise weak equivalence, and hence a trivial Reedy-fibration as required. Therefore, the asserted Bousfield localisation with respect to realisation weak equivalences exists.

For Segal-fibrant objects $X$, $Y$, the vertical arrows of diagram (\ref{Segal}) are trivial fibrations so that the underlying-space functor $U:\Top_{red}^{\Dop}\to\Top_*$ is homotopy-right-conservative; moreover, the canonical map $\phi:U(X)\to\Omega|X|_\Delta$ composed with the weak equivalence $\Omega|X|_\Delta\to|\Omega X|_\Delta$ coincides with the canonical map $i_0:(\Omega X)([0])\to|\Omega X|_\Delta$ under the identification $U(X)=X([1])=(\Omega X)([0])$. For Segal-fibrant $X$, the simplicial space $\Omega X$ is equifibered-fibrant, and therefore \emph{homotopically constant} (all simplicial operators act as weak equivalences).  If $X$ is also Reedy-cofibrant, then $\Omega X$ is pointwise cofibrant, which suffices to show that $i_0$, and hence $\phi$, is a weak equivalence. By Lemma \ref{nfold}, this proves that the localised model category of reduced simplicial spaces is a model for $1$-fold loop spaces.\end{proof}

\begin{rmk}\label{Segalfibrant2}Rezk \cite{R} and Bergner \cite{Ber} study related model structures on simplicial spaces. In both articles, only condition (a) of a Segal-fibration is taken into consideration; this leads to a different localisation of the pointwise model structure with a class of local equivalences which is strictly smaller than the class of realisation weak equivalences; the importance of this intermediate class comes from its relationship with the so-called \emph{Dwyer-Kan equivalences} of simplicially enriched categories. Rezk's concept of a \emph{complete} Segal-fibration is complementary to condition (b) of a Segal-fibration, in the following sense: a Reedy-fibration, which satisfies condition (a) of a Segal-fibration, is equifibered if and only if it satisfies at once condition (b) of a Segal-fibration and Rezk's completeness condition, cf. \cite[6.6]{R}.

It is an interesting problem to determine which properties of the monoidal model category $\Top$ and the cosimplicial object $\Delta_\Top^-:\Delta\to\Top$ are responsible for the validity of Theorem \ref{Segalfibrant}. To be more precise, let $\EE$ be an \emph{internal model category} (i.e. a cartesian closed category which is a monoidal model category with cofibrant unit for the cartesian product, cf. \cite{Ho}) such that the functor $|\!-\!|_\Delta:\EE^{\Dop}\to\EE$ induced by the cosimplicial object $\Delta_\EE^-:\Delta\to\EE$ is \emph{left exact}. All the constructions of this section make then sense and it is natural to ask under which conditions the commutative triangle of Quillen adjunctions\begin{gather}\begin{diagram}[noPS,p=1mm]\EE_*&\pile{\lTo^{|\!-\!|_{\Delta}}\\\rTo_{\Omega_{Seg}}}&\EE_{red}^{\Dop}\\&\pile{\rdTo^{\Omega}\\\luTo_{\Sigma}}&\dTo^U\uTo_L\\&&\EE_*\end{diagram}\end{gather}is a model for the derived image of the loop functor in $\EE_*$.

The proof of Theorem \ref{Segalfibrant} is quite formal and uses only two further properties which are specific to $\Top$. The first property is Theorem \ref{equifibered} which can be condensed in the following \emph{realisation axiom} of Rezk-Schwede-Shipley \cite{RSS}: Reedy-fibrations in $\EE^{\Dop}$ which are at once equifibered and realisation weak equivalences should be trivial Reedy-fibrations. The second property is the fibrancy of all objects which allows us deal with homotopy cartesian squares and to rephrase the homotopy-coreflection property of the Segal loop functor. This second property can however be weakened without any change in the proof of Theorem \ref{Segalfibrant}, by requiring that $\EE$ is \emph{right proper} and that the localised model category $\EE_{red}^{\Dop}$ is a \emph{good} model for the derived image of the loop functor; the latter amounts to the following \emph{fibrancy axiom}: the realisation functor $|\!-\!|_\Delta:\EE_{red}^{\Dop}\to\EE_*$ should take any cofibrant Segal-fibrant reduced simplicial object of $\EE$ to a cofibrant-fibrant object of $\EE_*$.\end{rmk}

\subsection{Reduced $\Gamma$-spaces}\label{Gammaspaces}

Segal's category $\Gamma$ is the category of finite sets $\nn=\{1,\dots,n\}$ (where $\underline{0}$ denotes the empty set) with operators $\mm\to\nn$ given by ordered $m$-tuples of pairwise disjoint subsets of $\nn$. Composition is defined by\begin{diagram}(\kk&\rTo^{(M_1,\dots,M_k)}&\mm&\rTo^{(N_1,\dots,N_m)}&\nn)=(\kk&\rTo^{(\bigcup_{j_1\in M_1}N_{j_1},\dots,\bigcup_{j_k\in M_k}N_{j_k})}&\nn).\end{diagram}A \emph{$\Gamma$-space} is a presheaf $X:\Gamma^\op\to\Top$. For us, a $\Gamma$-space is always supposed to be \emph{reduced}, i.e. $X(\underline{0})=*$. The underlying-space of a $\Gamma$-space $X$ is $U(X)=X(\underline{1})$. We briefly recall Segal's construction \cite{Se} of a spectrum out of a $\Gamma$-space, and review some features of Bousfield and Friedlander's \cite{BF} stable (injective) model structure on $\Gamma$-spaces. Segal's delooping functor $\Phi$ is part of the following commutative triangle of adjunctions:

\begin{gather}\begin{diagram}[noPS,silent,p=1mm]\Spt&\pile{\lTo^\Phi\\\rTo_\Psi}&\Top^{\Gamma^\op}_{red}\\&\pile{\rdTo^{\Omega^\infty}\\\luTo_{\Sigma^\infty}}&\dTo^U\uTo_L\\&&\Top_*\end{diagram}\end{gather}

Each based space $(X,*)$ defines a covariant functor $X^-:\Gamma\to\Top_*:\nn\mapsto X^n$. Indeed, it follows from the definitions that $\nn=\underline{1}^{\times n}$ in $\Gamma$, which implies that $\Top_*$ is equivalent to the category of \emph{left exact} functors $\Gamma\to\Top$ via $X\mapsto X^-$.

Each $\Gamma$-space $A$ induces an endofunctor $\underline{A}:\Top_*\to\Top_*:X\mapsto A\otimes_\Gamma X^-$. This endofunctor can also be defined as the left Kan extension of $A:\Gamma^\op\to\Top_*$ along $\Gamma^\op\inc\Top_*$, where $\Gamma^\op$ embeds in $\Top_*$ as the full subcategory spanned by the finite discrete sets $\{*,1,\dots,n\}$; a based map $\phi:\{*,1,\dots,n\}\to\{*,1,\dots,m\}$ is the dual of the $\Gamma$-operator $(\phi^{-1}(1)-\{*\},\dots,\phi^{-1}(m)-\{*\}):\mm\to\nn$.

To each $\Gamma$-space $A$ is associated (functorially in $A$) a binatural transformation\begin{gather}\begin{diagram}\phi_{X,Y}^A:\underline{A}(X)\wedge Y&\rTo&\underline{A}(X\wedge Y).\end{diagram}\end{gather}The \emph{Segal spectrum} $\Phi(A)$ is then defined by $\Phi(A)_n=\underline{A}(S^n)$ with structural maps $\phi_{S^n,S^1}^A:\Phi(A)_n\wedge S^1\to\Phi(A)_{n+1}$. The maps $\phi_{X,Y}^A$ are deduced from the canonical natural transformation $(X^-)\wedge Y\to(X\wedge Y)^-$. The right adjoint $\Psi$ of $\Phi$ uses the sphere-spectrum $\Sph=\Sigma^\infty S^0$ and the fact that $\Spt$ is enriched in $\Top$. Indeed, for any spectrum $E$, the $\Gamma$-space $\Psi(E)$ is defined by $\Psi(E)(\nn)=\underline{\Spt}(\Sph^{\times n},E)$, in particular $U\Psi=\Omega^\infty$ as required. For sake of completeness, we mention that the endofunctor $\underline{LX}$ equals $X\wedge -$, whence $\Phi L=\Sigma^\infty$ as required.

The \emph{stable} model category of $\Gamma$-spaces, as defined by Bousfield and Friedlander \cite{BF}, is the left Bousfield localisation with respect to stable equivalences of a \emph{strict} model category of $\Gamma$-spaces with pointwise weak equivalences. Similarly, the stable model category of spectra is the left Bousfield localisation with respect to stable equivalences of a strict model category of spectra with levelwise weak equivalences. The adjoint pair $(\Phi,\Psi)$ is a Quillen adjunction for the strict and the stable model categories. In order to see this, the following description of the structural maps $\phi_{S^n,S^1}^A$, essentially due to Segal \cite{Se}, will be very useful. 

Any $\Gamma$-space $A$ and based space $X$ define a $\Gamma$-space $A_X$ by $A_X=Am\otimes_\Gamma X^-$ where $m:\Gamma\times\Gamma\to\Gamma$ takes $(\nn_1,\nn_2)$ to $\underline{n_1n_2}$. Explicitly, the value of the $\Gamma$-space $A_X$ at $\nn$ is given by the coend formula$$A_X(\nn)=Am(\nn,-)\otimes_\Gamma X^-.$$As usual we denote by $m_!$ the left adjoint of the inverse image functor $m^*$. Kuhn, cf. \cite[Lemma 3.8]{K}, observes that for based spaces $X,Y$, there is a canonical identification $m_!(X^-\times Y^-)=(X\wedge Y)^-$; therefore, we get by adjointness\begin{align*}\underline{A_X}(Y)&=(Am\otimes_\Gamma X^-)\otimes_\Gamma Y^-=(m^*A)\otimes_{\Gamma\times\Gamma}(X^-\times Y^-)\\&=A\otimes_\Gamma m_!(X^-\times Y^-)=A\otimes_\Gamma(X\wedge Y)^-\\&=\underline{A}(X\wedge Y).\end{align*}Segal calls the $\Gamma$-space $A_{S^1}$ the \emph{classifying $\Gamma$-space} of $A$. Inductively, the $n$-fold classifying $\Gamma$-space of $A$ may be identified with $A_{S^n}$, and the Segal spectrum $\Phi(A)$ gets the formula $\Phi(A)_n=\underline{A_{S^n}}(S^0)=A_{S^n}(\underline{1})=U(A_{S^n})$, see \cite[Definition 1.3]{Se}. 

Segal also defines a functor $\gamma:\Delta\to\Gamma$ which takes the ordinal $[m]$ to the set $\mm$, and the simplicial operator $\phi:[k]\to[m]$ to the $\Gamma$-operator $\gamma(\phi):\kk\to\mm$ where $\gamma(\phi)(i)=\{j\in\mm\,|\,\phi(i-1)<j\leq\phi(i)\}$ for $i=1,\dots,k$. This yields for each $\Gamma$-space $A$ the identifications (cf. Corollary \ref{homeo} below and \cite[Proposition 3.2]{Se})\begin{align}\label{gamma1}|\gamma^*A|_\Delta=(\gamma^*A)\otimes_\Delta\Delta^-=A\otimes_\Gamma\gamma_!\Delta^-=A\otimes_\Gamma (S^1)^-=\underline{A}(S^1).\end{align}Bousfield and Friedlander \cite{BF} define for each $\Gamma$-space $A$ a \emph{skeletal filtration} $sk_nA$; the functor $\gamma^*$ takes this filtration to the skeletal filtration of the simplicial space $\gamma^*A$. The $1$-skeleton of $A_X$ is easily identified with $L(\underline{A}(X))$. Therefore, we get\begin{gather}\label{structural}\begin{diagram}|\gamma^*(sk_1A_X&\rTo&A_X)|_\Delta=(\underline{A}(X)\wedge S^1&\rTo^{\phi_{X,S^1}^A}&\underline{A}(X\wedge S^1)).\end{diagram}\end{gather}Since the composite functor $|\gamma^*(-)|_\Delta$ preserves colimits and (trivial) cofibrations of the strict model category of $\Gamma$-spaces, it follows immediately from (\ref{structural}) that the delooping functor $\Phi:\Top_{red}^{\Gamma^\op}\to\Spt$ preserves cofibrations and trivial cofibrations of the strict model categories. Therefore, $(\Phi,\Psi)$ is a Quillen adjunction for the strict model categories; it will be a Quillen adjunction for the stable model categories if $\Psi$ preserves stable fibrations. 

By an argument of Dugger \cite[8.5.4]{Hir}, it suffices to check that $\Psi$ preserves stable fibrations between stably fibrant objects. By Remark \ref{local}, stable fibrations between stably fibrant objects coincide with fibrations, and the latter are preserved by $\Psi$; so it suffices to show that $\Psi$ preserves stably fibrant objects. For a stably fibrant spectrum $E$, $\Psi(E)$ is a fibrant $\Gamma$-space which satisfies Segal's conditions (a) and (b), since for all $n$, the canonical inclusion of the $n$-fold wedge $\Sph^{\vee n}$ into the $n$-fold cartesian product $\Sph^{\times n}$ is a stable equivalence between cofibrant spectra. Therefore, it suffices to show that fibrant $\Gamma$-spaces which satisfy Segal's conditions (a) and (b) are stably fibrant. This may be deduced from the fact that $\Phi$ takes cofibrant-fibrant $\Gamma$-spaces $A$ which satisfy Segal's conditions (a) and (b) to $\Omega$-spectra. The latter property in turn is a consequence of Theorem \ref{Segalfibrant}, Lemma \ref{nfold} and formula (\ref{gamma1}), using that for all $n$, $\gamma^*(A_{S^n})$ is a cofibrant Segal-fibrant reduced simplicial space; indeed, this yields a canonical sequence of weak equivalences$$\underline{A}(S^0)\eqv\Omega\underline{A}(S^1)\eqv\Omega^2\underline{A}(S^2)\eqv\cdots.$$
\begin{thm}The category of reduced $\Gamma$-spaces with the stable injective model structure of Bousfield and Friedlander is a model for infinite loop spaces.\end{thm}
\begin{proof}All hypotheses of Lemma \ref{infinite} are fulfilled: we just showed that $(\Phi,\Psi)$ is a Quillen adjunction, and that $\Phi$ takes cofibrant-stably fibrant $\Gamma$-spaces to $\Omega$-spectra. It is clear that $(L,U)$ also is a Quillen adjunction, and that $U$ is homotopy-right-conservative. Finally, we have equality $U=\Omega^\infty\Phi$.\end{proof}

\begin{rmk}\label{connectivespectra}The previous theorem is just a reformulation of one of the central results of Segal \cite{Se} and Bousfield-Friedlander \cite{BF}, namely that reduced $\Gamma$-spaces are models for connective spectra. Indeed, the homotopy-counits of the $\Phi$-$\Psi$-adjunction yield an explicit connective cover for each $\Omega$-spectrum. The homotopy category of reduced $\Gamma$-spaces may thus be identified with the coreflective subcategory of the stable homotopy category spanned by the connective spectra.\end{rmk}

\section{Wreath product over $\Delta$ and duality.}

In this central section, we define for each small category $\AA$, categorical wreath products $\Delta\wr\AA$ and $\Gamma\wr\AA$. The main object of study are the iterated wreath products $\Delta\wr\cdots\wr\Delta$ which will be denoted by $\Theta_n$. We use this notation, since we shall show below that the dual of $\Theta_n$ may be identified with Joyal's category $\DD_n$ of finite combinatorial $n$-disks, and Joyal \cite{Joy1} uses the notation $\Theta_n$ for the dual of $\DD_n$. Our point of view here is to consider $\Theta_n$ as an $n$-categorical analog of $\Delta$, and to exploit the fact that the wreath product construction allows us to argue by induction on $n$, building on classical properties of the simplex category $\Delta=\Theta_1$. In particular, we shall establish the following four properties:\begin{enumerate}\item $\Theta_n$ is a dense subcategory of the category of (small strict) $n$-categories;\item the presheaf topos on $\Theta_n$ is a classifying topos for combinatorial $n$-disks;\item the presheaf topos on $\Theta_n$ has a left exact realisation functor with values in the category of $\CW$-complexes;\item there is an assembly functor $\gamma_n:\Theta_n\to\Gamma$ whose induced geometric morphism classifies the generic $n$-sphere (i.e. the quotient of the generic combinatorial $n$-disk by its boundary).\end{enumerate}

Since all idempotents of $\Theta_n$ split, property (2) implies\vspace{1ex}\begin{enumerate}\item[(2$'$)]the dual of $\Theta_n$ is isomorphic to the category of \emph{finite} combinatorial $n$-disks.\end{enumerate}The duality between $\Theta_n$ and $\DD_n$, as expressed by (2) and (2$'$), has an analog for Segal's category $\Gamma$: the presheaf topos on $\Gamma$ is a classifying topos for pointed objects; this implies that the dual of $\Gamma$ is the category of \emph{finite pointed sets}, cf. Section \ref{Gammaspaces}, and shows that (4) characterises the assembly functor $\gamma_n$.

For any small category $\AA$, the presheaf topos $\Sets^{\AA^\op}$ will be denoted by $\widehat{A}$; for any object $A$ of $\AA$, the representable presheaf $\AA(-,A)$ will be denoted by $\AA[A]$.

\begin{dfn}\label{def}The wreath product $\Delta\wr\AA$ (resp. $\Gamma\wr\AA$) is defined to be the category\vspace{1ex}
 
$\bullet$ with objects the $m$-tuples $(A_1,\dots,A_m)$ of objects of $\AA$, for varying $m\geq 0$;

$\bullet$ with operators $(\phi;\phi_1,\dots,\phi_m):(A_1,\dots,A_m)\to(B_1,\dots,B_n)$ all $(m+1)$-tuples  consisting of a $\Delta$-operator $\phi:[m]\to[n]$ (resp. $\Gamma$-operator $\phi:\mm\to\nn$) and an $m$-tuple $(\phi_1,\dots,\phi_m)$ of morphisms in $\widehat{\AA}$ of the form $$\phi_i:\AA[A_i]\to\AA[B_{\phi(i-1)+1}]\times\AA[B_{\phi(i-1)+2}]\times\cdots\times\AA[B_{\phi(i)}]$$ $$(\textrm{resp. }\phi_i:\AA[A_i]\to\prod_{k\in\phi(i)}\AA[B_k]).$$\end{dfn}
\noindent If the set $\{k\,|\,\phi(i-1)<k\leq\phi(i)\}$ (resp. $\phi(i)$) is empty, the target of $\phi_i$ above is the terminal presheaf. If $[m]=[0]$ (resp. $\mm=\underline{0}$), operators $(\phi;):()\to(B_1,\dots,B_n)$ in $\Delta\wr\AA$ (resp. $\Gamma\wr\AA$) correspond bijectively to operators $\phi:[0]\to[n]$ (resp. $\underline{0}\to\nn$). Therefore, the $0$-tuple in $\Delta\wr\AA$ (resp. $\Gamma\wr\AA$) is the terminal object of $\Delta\wr\AA$ (resp. null-object of $\Gamma\wr\AA$), since $[0]$ is the terminal object of $\Delta$ (resp. $\underline{0}\,$ is the null-object of $\Gamma$). Recall that by definition a \emph{null-object} of a category is an object which is simultaneously initial and terminal.

Composition of operators is obvious: indeed, the use of the Yoneda-embedding $\AA\inc\widehat{\AA}$ above may be avoided. A more elementary (but more cumbersome) definition replaces the morphism $\phi_i:\AA[A_i]\to\AA[B_{\phi(i-1)+1}]\times\AA[B_{\phi(i-1)+2}]\times\cdots\times\AA[B_{\phi(i)}]$ by its components $\phi_i^k:A_i\to B_k$ in $\AA$; composition in the wreath product is then directly induced by composition in $\AA$. 

If $\AA$ has a null-object, the category $\AA\,$ may be considered as a pointed object in $\Cat$ with base point $\star\to\AA$ given by the null-object of $\AA$. Therefore, there is a canonical functor $\AA^-:\Gamma\to\Cat$ taking $\nn$ to $\AA^{\times n}$, cf. Section \ref{Gammaspaces}. The \emph{Grothendieck-construction} $\int_\Gamma\AA^-$ is then isomorphic to the wreath product $\Gamma\wr\AA$. Similarly, for the Segal functor $\gamma:\Delta\to\Gamma$ composed with $\AA^-:\Gamma\to\Cat$, the Grothendieck-construction $\int_\Delta\AA^-\circ\gamma$ is isomorphic to the wreath product $\Delta\wr\AA$. However, if $\AA$ has no null-object, it seems unlikely that the wreath products $\Gamma\wr\AA$ and $\Delta\wr\AA$ may be described by such Grothendieck-constructions.

The Segal functor $\gamma:\Delta\to\Gamma$, cf. Section \ref{Gammaspaces}, induces an obvious functor $\gamma\wr\AA:\Delta\wr\AA\to\Gamma\wr\AA$ which is the identity on objects, and sends an operator $(\phi;\phi_1,\dots,\phi_m)$ of $\Delta\wr\AA$ to the operator $(\gamma(\phi);\phi_1,\dots,\phi_m)$ of $\Gamma\wr\AA$. We shall see that the two wreath products enjoy quite different ``universal properties''; the wreath product over $\Delta$ is a kind of non-symmetric version of the wreath product over $\Gamma$.

Just for the purpose of the next proposition we introduce the following terminology: a category $\CC$ will be called \emph{semi-additive} whenever $\CC$ has finite products and a null-object. A finite product-preserving functor between semi-additive categories will be called semi-additive.

\begin{lma}\label{semi-additive}The wreath product $\Gamma\wr\AA$ is the free semi-additive category on $\AA$. In particular, there is a canonical assembly functor $\alpha:\Gamma\wr\Gamma\to\Gamma$ which takes $(\nn_1,\cdots,\nn_k)$ to $\underline{n_1+\cdots+n_k}$.\end{lma}
\begin{proof}By definition, $\underline{0}$ is a null-object of $\Gamma$ so that the $0$-tuple $()$ is a null-object for $\Gamma\wr\AA$. Concatenation of tuples defines a cartesian product for $\Gamma\wr\AA$ so that $\Gamma\wr\AA$ is semi-additive; moreover, $\AA$ embeds in $\Gamma\wr\AA$ via $A\mapsto (A)$, and the universal property of $\Gamma\wr\AA$ follows. The assembly functor $\alpha$ is induced by the universal property of $\Gamma\wr\Gamma$ applied to the identity functor $\Gamma\to\Gamma\wr *$.\end{proof}

In particular, each functor $F:\AA\to\BB$ induces a canonical semi-additive functor $\Gamma\wr F:\Gamma\wr\AA\to\Gamma\wr\BB$ taking $(A)$ to $(F(A))$. This yields a composite functor\begin{gather*}\begin{diagram}\gamma\wr F:\Delta\wr\AA&\rTo^{\gamma\wr\AA}&\Gamma\wr\AA&\rTo^{\Gamma\wr F}&\Gamma\wr\BB.\end{diagram}\end{gather*}

\begin{dfn}\label{nerve}Let $\,\Theta_1=\Delta$ and $\gamma_1=\gamma$, cf. \ref{Gammaspaces}. Define by induction on $n$, $$\Theta_n=\Delta\wr\Theta_{n-1},\textrm{ and}$$\begin{diagram}\gamma_n:\Theta_n=\Delta\wr\Theta_{n-1}&\rTo^{\gamma\wr\gamma_{n-1}}&\Gamma\wr\Gamma&\rTo^\alpha&\Gamma.\end{diagram}\end{dfn}

\begin{rmk}The preceding definition of $\Theta_n$ applies the wreath product in a specific order, namely $\Theta_n=\Delta\wr(\Delta\wr(\cdots\wr(\Delta\wr\Delta)\cdots))$. We are grateful to J.-L. Loday for pointing out to us the following generalisation of the wreath product which allows us to formulate an \emph{associativity property} of the wreath product.

It is straightforward to check that all we need to carry out the definition of $\BB\wr\AA$ for arbitrary small categories $\BB$ and $\AA$ is the existence of a specific functor $\BB\to\Gamma$. If $\AA$ also comes equipped with a functor to $\Gamma$, the wreath product $\BB\wr\AA$ is again equipped with a functor to $\Gamma$, exactly as above; in other words, the wreath product is actually a bifunctor $-\wr-:\Cat/\Gamma\times\Cat/\Gamma\lra\Cat/\Gamma.$ It can be checked that this bifunctor defines a monoidal structure on $\Cat/\Gamma$ with unit the functor $\star\to\Gamma$ taking the unique object of the one-point category to $\underline{1}$. In particular, if we consider $\Theta_n$ as equipped with the above defined assembly functor $\gamma_n:\Theta_n\to\Gamma$, then we get canonical isomorphisms of categories $\Theta_m\wr\Theta_n\cong\Theta_{m+n}$. We shall however not need these associativity isomorphisms in this article.\end{rmk}

For the next proposition, recall that a small subcategory $\AA$ of a category $\VV$ is called \emph{dense} (or \emph{adequate} in the original terminology of Isbell) if the associated nerve functor $\NN_\AA:\VV\to\widehat{\AA}:V\mapsto\VV(-,V)$ is fully faithful.

\begin{prp}\label{wreath}For any full (small) subcategory $\AA$ of a cocomplete cartesian closed category $\VV$, the wreath product $\Delta\wr\AA$ embeds in $\VV$-$\Cat$ as the full subcategory spanned by the free $\VV$-categories on finite linear $\AA$-graphs. 

If $\AA$ is dense in $\VV$ then $\Delta\wr\AA$ is dense in $\VV$-$\Cat$.\end{prp}

\begin{proof}Under the given hypotheses, the forgetful functor $U_\VV:\VV$-$\Cat\to\VV$-$\Grph$ from $\VV$-categories to $\VV$-graphs has a left adjoint $F_\VV:\VV$-$\Grph\to\VV$-$\Cat$. A linear $\VV$-graph of length $m$ is just an $m$-tuple $(V_1,\dots,V_m)$ of objects of $\VV$, so that there is a canonical correspondence between the objects of $\Delta\wr\AA$ and finite linear $\AA$-graphs. The free $\VV$-category on $(A_1,\dots,A_m)$ admits then the following explicit description:$$F_\VV(A_1,\dots,A_m)(i,j)=A_{i+1}\times\cdots\times A_j$$where $0\leq i<j\leq m$. The endomorphism-objects are terminal in $\VV$, all other hom-objects are initial in $\VV$. Composition is induced by cartesian product in $\VV$. It follows that the set of $\VV$-functors $\FF_\VV(A_1,\dots,A_m)\to\FF_\VV(B_1,\dots,B_n)$ may be identified with the set of operators $(A_1,\dots,A_m)\to(B_1,\dots,B_n)$ in $\Delta\wr\AA$, and that this correspondence respects the composition laws.

For the second statement, we have to show that if the nerve functor $\NN_\AA:\VV\to\widehat{A}$ is fully faithful, then the nerve functor $\NN_{\Delta\wr\AA}:\VV$-$\Cat\to\widehat{\Delta\wr\AA}$ is also fully faithful. For this, we factor the latter as follows:\begin{diagram}[noPS,small,silent]\VV\!-\!\Cat&\rTo^{\NN_{\Delta\wr\AA}}&\widehat{\Delta\wr\AA}\\\dTo^{\NN_\AA\!-\!\Cat}&\ruTo_{\NN'_{\Delta\wr\AA}}&\\\widehat{\AA}\!-\!\Cat&&\end{diagram}where the vertical functor is a ``change-of-enrichment'' functor and $\NN'_{\Delta\wr\AA}$ is derived from the construction above applied to the Yoneda-embedding $\AA\inc\widehat{\AA}$. The vertical functor is fully faithful since by hypothesis $\NN_\AA$ is fully faithful. Therefore, it suffices to show that $\NN'_{\Delta\wr\AA}$ is fully faithful: for any $\widehat{\AA}$-category $X$, the value of $\NN'_{\Delta\wr\AA}(X)$ at $(A_1,\dots,A_m)$ is given by$$\coprod_{(x_0,\dots,x_m)\in X_0^{m+1}}X(x_0,x_1)(A_1)\times X(x_1,x_2)(A_2)\times\cdots\times X(x_{m-1},x_m)(A_m)$$where $X_0$ denotes the set of objects of $X$.

Setting $m=1$ shows that $\NN'_{\Delta\wr\AA}$ is faithful. For any $\widehat{\Delta\wr\AA}$-morphism $f:\NN'_{\Delta\wr\AA}(X)\to\NN'_{\Delta\wr\AA}(Y)$, the compatibility of $f$ with $\Delta\wr\AA$-operators of type $(id_{[1]};\phi_1):(A)\to(B)$ yields an $\widehat{\AA}$-morphism $X(x,y)\to Y(fx,fy)$ for each pair of objects $(x,y)\in X_0^2$, the compatibility with $\Delta\wr\AA$-operators of type $(\phi;\phi_1,\phi_2):(A_3)\to(A_1,A_2)$ shows that $f$ is the nerve of a $\widehat{A}$-functor, thus $\NN'_{\Delta\wr\AA}$ is full.\end{proof}

Proposition \ref{wreath} inductively identifies $\Theta_n$ with a dense subcategory of the category $\nCat$ of (small strict) $n$-categories. In order to make this embedding explicit, we need Batanin's star-construction which establishes a one-to-one correspondence between \emph{level-trees} $T$ of height $\leq n$ and certain \emph{$n$-graphs} $T_*$, cf. \cite[p. 61]{Ba1}, \cite[1.2]{mi}. 

\subsection{Level-trees and hypergraphs}\label{level-tree}A \emph{level-tree} $T$, cf. \cite[1.1-4]{mi}, is a contractible planar graph with a distinguished root-vertex, and (if existent) a distinguished ``left-most'' root-edge; the \emph{height} of a vertex $v$ of $T$ is the length of the edge-path between $v$ and the root-vertex; the \emph{height} of $T$ is the maximal occuring height of its vertices; the edges are rootwards oriented so that for each vertex of $T$ the incoming edges are totally ordered from left to right according to the planarity of $T$ and to the distinguished left-most root-edge. 

An \emph{$n$-graph} $X$, cf. \cite[0.2]{mi}, is a graded set $(X_k)_{0\leq k\leq n}$ together with source/target maps $s,t:X_k\dto X_{k-1}$ satisfying the usual relations $ss=st$ and $ts=tt$. The category of $n$-graphs is thus a presheaf category which we shall denote by $\nGrph$. 

Street \cite{St2} has shown that the $n$-graphs $X$ arising as $T_*$ for some level-tree $T$ are those for which the partial order on $\bigsqcup_{k=0}^nX_k$, generated by $s(x)\leq x\leq t(x)$, is a \emph{total order}. In particular, the representable $n$-graphs are the $n$-graphs $\bar{k}_*$ for linear level-trees $\bar{k}$ of height $k\leq n$. In general, $T_*$ is a canonical ``$T$-shaped'' colimit of representable $n$-graphs, one for each ``sector'' of $T$, cf. \cite[pg. 61]{Ba1}, \cite[1.2]{mi}.

The objects of the iterated wreath product $\Theta_n=\Delta\wr\Theta_{n-1}$ are ``bunches'' of objects of $\Theta_{n-1}$ or, equivalently, level-trees of height $\leq 1$ the edges of which support objects of $\Theta_{n-1}$. Assuming inductively that the objects of $\Theta_{n-1}$ are finite level-trees of height $\leq n-1$, this establishes a one-to-one correspondence between the objects of $\Theta_n$ and finite level-trees of height $\leq n$. We shall henceforth identify the objects of $\Theta_n$ with finite level-trees of height $\leq n$ under this correspondence.

The obvious forgetful functor $U_n:\nCat\to\nGrph$ has a left adjoint $F_n:\nGrph\to\nCat$ which can be described explicitly, cf. \cite{Ba1}, \cite[1.8, 1.12]{mi}.

\begin{thm}\label{dense}The iterated wreath product $\,\Theta_n$ embeds in $\nCat$ as the full subcategory spanned by the free $n$-categories $F_n(T_*)$ on $T_*$, where $T$ runs through the set of finite level-trees of height $\leq n$. In particular, $\Theta_n$ is dense in $\nCat$.\end{thm}

\begin{proof}We proceed by induction on $n$. For $n=1$, the statement is well known: level-trees of height $\leq 1$ are finite (possibly empty) corolla $T^m$ which induce the linear $1$-graphs $T^m_*$ of length $m$. The simplex-category $\Delta=\Theta_1$ is the full subcategory of $\Cat$ spanned by the finite ordinals $[m]=F_1(T^m_*)$. The density of $\Delta$ in $\Cat$ is (well known and) the special case $\AA=*$ of Proposition \ref{wreath}. 

Consider now for $\VV=\nCat$ and $\WW=\nGrph$ the following factorisation of the adjunction $U_{n+1}:\nplusCat\rlto\nplusGrph:F_{n+1}$:\begin{diagram}\nplusCat=\VV\!-\!\Cat&\pile{\rTo^{U_\VV}\\\lTo_{F_\VV}}&\VV\!-\!\Grph&\pile{\rTo^{U_n\!-\!\Grph}\\\lTo_{F_n\!-\!\Grph}}&\WW\!-\!\Grph=\nplusGrph.\end{diagram}The inductive step follows then from Proposition \ref{wreath} and the fact that under the identification $\nplusGrph=\WW\!-\!\Grph$, Batanin's star-construction $T_*$ of a level-tree $T$ of height $\leq n+1$ and root-valence $m$ corresponds to the $1$-graph $T^m_*$ labelled by the $m$-tuple of $n$-graphs $(S^1_*,\dots,S^m_*)$, where $T$ is obtained by grafting the $m$-tuple $(S^1,\dots,S^m)$ onto the corolla $T^m$.\end{proof}

This proves property (1) of the introduction to this section. For property (2), recall (cf. \cite[III.4]{Bor}) that a presheaf topos $\widehat{\AA}$ is called a \emph{classifying topos} for a certain structure $\mathcal{T}$ if there is an equivalence of categories between the category of flat functors $\AA\to\Sets$ and the category of $\mathcal{T}$-structured sets. A functor $F:\AA\to\Sets$ is called \emph{flat} if its left Kan extension $(-)\otimes_\AA F:\widehat{\AA}\to\Sets$ is left exact and therefore the inverse-image-part of a geometric morphism $\Sets\rlto\widehat{\AA}$. Precomposing the inverse-image-part with the Yoneda-embedding $\AA\inc\widehat{\AA}$ gives back the flat functor $\AA\to\Sets$ that represents the $\mathcal{T}$-structured set. Therefore, the Yoneda-embedding $\AA\inc\widehat{\AA}$ is also called a \emph{generic model} for $\mathcal{T}$-structured objects. Property (2) thus states that there is an equivalence of categories between the category of flat functors $\Theta_n\to\Sets$ and the category of combinatorial $n$-disks in Joyal's sense \cite{Joy1}.

We shall prove property (2) by induction on $n$. Again, the case $n=1$ is well known, since a combinatorial $1$-disk is defined to be a linearly ordered set $(X,\leq_X)$ with \emph{distinct} minimal and maximal elements $x_0$ and $x_1$, and any such quadruple $(X,\leq_X,x_0,x_1)$ determines, and is determined by, a flat functor $\underline{X}:\Delta\to\Sets$ such that $(X,x_0,x_1)=\underline{X}([0]\dto[1])$. The linear order $\leq_X$ may be recovered from $\underline{X}$ as the relation on $X$ induced by the degeneracy pair $(s^*_0,s^*_1):\underline{X}([2])\to\underline{X}([1])\times\underline{X}([1])$, cf. MacLane and Moerdijk \cite[VIII.8]{MM}.

\begin{dfn}\label{diagonal}For each small category $\AA$, define $\delta_\AA:\Delta\times\AA\to\Delta\wr\AA$ by\begin{align*}\delta_\AA([n],A)&=\overbrace{(A,\dots,A)}^\text{n times}\\\delta_\AA(\phi:[m]\to[n],f:A\to B)&=(\phi;f_1,\dots,f_m)\end{align*} where $f_i:\AA[A]\to\AA[B_{\phi(i-1)+1}]\times\cdots\times\AA[B_{\phi(i)}]$ is $\AA[f]:\AA[A]\to\AA[B]$ on each factor.

Define inductively $\delta_1:\Delta\overset{id}{\lra}\Theta_1$ and for $n>1$,\begin{diagram}[noPS,small]\delta_n&:\Delta^{\times n}&\rTo^{\delta_{\Delta^{\times n-1}}}&\Delta\wr\Delta^{\times n-1}&\rTo^{\Delta\wr\delta_{n-1}}&\Delta\wr\Theta_{n-1}=\Theta_n.\end{diagram}\end{dfn}

For the inductive step in the proof of property (2) we need a preliminary result on realisation functors for the presheaf topos $\widehat{\Delta\wr\AA}$. By a \emph{realisation functor}, we mean any colimit-preserving \emph{left-exact} functor, e.g. the inverse-image-part of a geometric morphism of toposes.

\begin{prp}\label{realisation}Let $\EE$ be a cocomplete cartesian closed category, and assume that $\widehat{\Delta}$ and $\widehat{\AA}$ have $\EE$-valued realisation functors $|\!-\!|_\Delta$ and $|\!-\!|_\AA$. Then $\widehat{\Delta\wr\AA}$ has an $\EE$-valued realisation functor $|\!-\!|_{\Delta\wr\AA}=|\!-\!|_{\Delta\times\AA}\circ(\delta_\AA)^*$ taking $(\Delta\wr\AA)[(A_1,\dots,A_m)]$ to the \emph{simplicial suspension} of the $m$-tuple $(|\AA[A_1]|_\AA,\dots,|\AA[A_m]|_\AA)$.\end{prp}

\begin{proof}The presheaf topos $\widehat{\Delta\times\AA}$ has an obvious realisation functor $|X\times Y|_{\Delta\times\AA}=|X|_\Delta\times|Y|_\AA$. Therefore, since $(\delta_\AA)^*$ is left and right adjoint, the composite functor $|\!-\!|_{\Delta\times\AA}\circ(\delta_\AA)^*$ defines a realisation functor $|\!-\!|_{\Delta\wr\AA}$ for $\widehat{\Delta\wr\AA}$.

It remains to be shown that on representable objects, this realisation functor is a simplicial suspension functor. In order to give a precise definition of the latter, observe first that the given realisation functor $|\!-\!|_\Delta:\widehat{\Delta}\to\EE$ factors through the category $\EE^{\Dop}$ of simplicial $\EE$-objects via a functor $(-)_\EE:\Sets^{\Dop}\to\EE^{\Dop}$, obtained by composition with the obvious functor $\Sets\to\EE:X\mapsto\bigsqcup_X\star_\EE$, where $\star_\EE$ denotes the terminal object of $\EE$. In other words, if we denote by $\Delta^-_\EE:\Delta\to\EE$ the composite functor $\Delta\inc\widehat{\Delta}\to\EE$, then we have canonical identifications $\Delta^n_\EE=|\Delta[n]_\EE|_\Delta$. Next, observe that for each object $E$ of $\EE$, there is an essentially unique map of simplicial $\EE$-objects $\Delta[E]\to\Delta[1]_\EE$ such that the fibers over the extremities $\Delta[0]_\EE\dto\Delta[1]_\EE$ are $*_\EE$ and all other fibers are $E$. The \emph{unreduced suspension} $SE$ of $E$ is by definition the realisation of the simplicial $\EE$-object $\Delta[E]$; the unreduced suspension is thus fibered over $\Delta^1_\EE$. More generally, the \emph{simplicial suspension} $S(E_1,\dots,E_m)$ of an $m$-tuple $(E_1,\dots,E_m)$ is the realisation of the simplicial $\EE$-object $\Delta[E_1,\dots,E_m]$ defined by the following pullback in $\EE^{\Dop}$:\begin{gather}\label{suspension}\begin{diagram}[noPS,small]\Delta[E_1,\dots,E_m]&\rTo&\Delta[E_1]\times\cdots\times\Delta[E_m]\\\dTo&&\dTo\\\Delta[m]_\EE&\rTo&\Delta[1]_\EE\times\cdots\times\Delta[1]_\EE,\end{diagram}\end{gather}where the lower horizontal map is induced by the canonical inclusion of $\Delta[m]$ into the simplicial $m$-cube $\Delta[1]^m$ (ordering the $m$ factors from left to right). In particular, the simplicial suspension $S(E_1,\dots,E_m)$ is fibered over $\Delta^m_\EE$.

We want to show that $|(\Delta\wr\AA)[(A_1,\dots,A_m)]|_{\Delta\wr\AA}=S(|\AA[A_1]|_\AA,\dots,|\AA[A_m]|_\AA).$ Any object of $\widehat{\Delta\times\AA}$ can be realised in two steps, by realising first with respect to the $\AA$-coordinate, and then with respect to the $\Delta$-coordinate; in particular, the realisation functor $|\!-\!|_{\Delta\wr\AA}:\widehat{\Delta\wr\AA}\to\EE$ factors through the category of simplicial $\EE$-objects via the composite  functor $r_\AA:\widehat{\Delta\wr\AA}\overset{(\delta_\AA)^*}{\lra}\widehat{\Delta\times\AA}\overset{|\!-\!|_\AA}{\lra}\EE^{\Dop}$. Therefore, it remains to be shown that the simplicial $\EE$-object $\Delta[|\AA[A_1]|_\AA,\dots,|\AA[A_m]|_\AA]$ may be identified with $r_\AA(\Delta\wr\AA)[(A_1,\dots,A_m)]$. We show more precisely that the pullback square (\ref{suspension}) is the image under $r_A$ of the following commutative square\begin{gather}\label{suspension2}\begin{diagram}[noPS,small](\Delta\wr\AA)[(A_1,\dots,A_m)]&\rTo^\beta&(\Delta\wr\AA)[(A_1)]\times\cdots\times(\Delta\wr\AA)[(A_m)]\\\dTo^\alpha&&\dTo\\\pi^*\pi_!(\Delta\wr\AA)[(A_1,\dots,A_m)]&\rTo&\pi^*\pi_!(\Delta\wr\AA)[(A_1)]\times\cdots\times\pi^*\pi_!(\Delta\wr\AA)[(A_m)]\end{diagram}\end{gather}where $\pi:\Delta\wr\AA\to\Delta$ is the canonical projection, $E_k=|\AA[A_k]|_\AA$ for $1\leq k\leq m$ and the vertical maps are induced by the units of the $(\pi_!,\pi^*)$-adjunction. The horizontal maps of (\ref{suspension2}) are induced by the $m$ surjective $\Delta$-operators $\phi:[m]\to[1]$. Since the left adjoint $\pi_!$ preserves representable presheaves, the lower horizontal map of (\ref{suspension2}) is the canonical inclusion $\pi^*\Delta[m]\inc\pi^*\Delta[1]\times\cdots\times\pi^*\Delta[1]$ whose image under $r_\AA$ is precisely the lower horizontal map of (\ref{suspension}). It follows from the definitions that for any object $A$ of $\AA$, the image under $r_\AA$ of the representable presheaf $(\Delta\wr\AA)[(A)]$ is the simplicial $\EE$-object $\Delta[|\AA[A]|_\AA]$ defined above; in particular, the image of the right vertical map of (\ref{suspension2}) is the right vertical map of (\ref{suspension}). Therefore, since $r_\AA$ is left exact, it finally remains to be shown that (\ref{suspension2}) is a pullback square in $\widehat{\Delta\wr\AA}$. But this follows from the definition of the $(\Delta\wr\AA)$-operators: indeed, for any $(\phi;\phi_1,\dots,\phi_n):(B_1,\dots,B_n)\to(A_1,\dots,A_m)$, considered as an element of $(\Delta\wr\AA)[(A_1,\dots,A_m)]$, the image under $\alpha$ determines $\phi$, the image under $\beta$ determines $(\phi_1,\dots,\phi_n)$, and the commutativity of (\ref{suspension2}) amounts to the property that $\phi_i$ represents a map $\AA[B_i]\to\AA[A_{\phi(i-1)+1}]\times\cdots\times\AA[A_{\phi(i)}]$ in $\widehat{\AA}$.\end{proof}


\begin{thm}\label{disk}The presheaf topos on the iterated wreath peoduct $\,\Theta_n$ is a classifying topos for combinatorial $n$-disks. In particular, the dual of $\,\Theta_n$ is isomorphic to Joyal's category $\DD_n$ of finite combinatorial $n$-disks.\end{thm}

\begin{proof}As mentioned above, we argue by induction on $n$. Recall from \cite{Joy1} that a \emph{combinatorial $n$-disk} consists of a sequence of sets and mappings\begin{diagram}[noPS,small,silent]D_n&\pile{\lTo^{s_n}\\\rTo^{i_n}\\\lTo^{t_n}}&D_{n-1}&\pile{\lTo^{s_{n-1}}\\\rTo^{i_{n-1}}\\\lTo^{t_{n-1}}}&D_{n-2}&\cdots&D_2&\pile{\lTo^{s_2}\\\rTo^{i_2}\\\lTo^{t_2}}&D_1&\pile{\lTo^{s_1}\\\rTo^{i_1}\\\lTo^{t_1}}&D_0\end{diagram}such that $D_0$ is singleton, $i_ks_k=i_kt_k=id_{D_{k-1}}$, $t_{k+1}t_k=s_{k+1}t_k$, $t_{k+1}s_k=s_{k+1}s_k$ for all $k$, and such that for each $x\in D_{k-1}$, the fiber $i_k^{-1}(x)$ is linearly ordered with minimum $s_k(x)$ and maximum $t_k(x)$, and finally the equalizer of $(s_k,t_k)$ is the empty set for $k=1$, resp. the union $s_{k-1}(D_{k-2})\cup t_{k-1}(D_{k-2})$ for $k=2,\dots,n$. 

Therefore, if $\bar{\GG}_n$ denotes the full subcategory of $\Theta_n$ spanned by the linear level-trees, cf. \cite[0.2/2.1]{mi}, a combinatorial $n$-disk amounts to covariant functor $D^n:\bar{\GG}_n\to\Sets$ fullfilling certain order and exactness relations. Maps of combinatorial $n$-disks are natural transformations preserving the linear orderings of the fibers. We have to show that any such combinatorial $n$-disk $\bar{\GG}_n\to\Sets$ extends in an essentially unique way to a flat functor $\Theta_n\to\Sets$ and, conversely, any flat functor $\Theta_n\to\Sets$ restricts to a combinatorial $n$-disk on $\bar{\GG}_n$. 

The main idea for the inductive step is to cut the $n$-disk $D^n$ into a $1$-disk $D^1:\bar{\GG}_1\to\Sets$ and a ``desuspended'' $(n-1)$-disk $D^{n-1}:\bar{\GG}_{n-1}\to\Sets$. More precisely, the combinatorial $1$-disk $D^1$ is obtained as the restriction of $D^n$ along the canonical inclusion $\bar{\GG}_1\inc\bar{\GG}_n$; in particular, since simplicial sets form a classifying topos for linearly ordered sets, this $1$-disk extends to a flat functor $\Delta\to\Sets$ and (by left Kan extension) to a realisation functor $\widehat{\Delta}\to\Sets$. Therefore, we get by Proposition \ref{realisation} (with $\AA=\star$ and $\EE=\Sets$) an \emph{unreduced suspension functor} $S:\Sets\to\Sets$. Due to the exactness property of an $n$-disk, there is then an essentially unique $(n-1)$-disk $D^{n-1}$ such that the value of $D^n$ at a linear level-tree of height $k>0$ equals the value of $SD^{n-1}$ at a linear level-tree of height $k-1$. Now, apply the induction hypothesis and extend $D^{n-1}$ to a flat functor on $\Theta_{n-1}$ and (by left Kan extension) to a realisation functor for $\widehat{\Theta}_{n-1}$. Hence, according to Proposition \ref{realisation}, we get a realisation functor for $\widehat{\Theta}_n$ which on representable objects is a simplicial suspension functor. Restricting back along the Yoneda-embedding $\Theta_n\inc\widehat{\Theta}_n$ yields the required essentially unique flat extension of $D^n$. For the converse direction, a similar inductive argument applies.

For the duality, we use that a combinatorial $n$-disk $\bar{\GG}_n\to\Sets$ is finite if and only if its flat extension $\Theta_n\to\Sets$ preserves filtered colimits. Moreover, idempotents in $\Theta_n$ split, cf. \cite[Remark 2.5]{mi} and Remark \ref{factor}. As a consequence, a filtered colimit preserving flat functor $\Theta_n\to\Sets$ is \emph{corepresentable}, cf. \cite[I:6.5.6/6.7.6]{Bor}); thus, the contravariant Yoneda-embedding $\Theta_n^\op\inc\Sets^{\Theta_n}:T\mapsto\Theta(T,-)$ identifies the dual of $\,\Theta_n$ with Joyal's category $\DD_n$ of finite combinatorial $n$-disks.\end{proof}

\begin{cor}\label{flat}The presheaf topos on $\,\Theta_n$ admits a topological realisation functor which takes the representable presheaf $\,\Theta_n[T]$ to a $\CW$-ball of dimension equal to the number of edges of the level-tree $T$.\end{cor}

\begin{proof}This follows either inductively from Proposition \ref{realisation}, using the classical realisation functor for simplicial sets, or directly from Theorem \ref{disk}, using that $\Top$ is topological over $\Sets$ (see \cite[II:7.3.2]{Bor}), and that the underlying set of the unit $n$-ball $B^n$ in $\RR^n$ has a combinatorial $n$-disk structure, cf. \cite{Joy1}, \cite[2.1]{mi}:\begin{diagram}[noPS,small,silent]B^n&\pile{\lTo^{s_n}\\\rTo^{i_n}\\\lTo^{t_n}}&B^{n-1}&\pile{\lTo^{s_{n-1}}\\\rTo^{i_{n-1}}\\\lTo^{t_{n-1}}}&B^{n-2}&\cdots&B^2&\pile{\lTo^{s_2}\\\rTo^{i_2}\\\lTo^{t_2}}&B^1&\pile{\lTo^{s_1}\\\rTo^{i_1}\\\lTo^{t_1}}&B^0\end{diagram}where $s_k$ (resp. $t_k$) maps $B^{k-1}$ to the lower (resp. upper) hemisphere of $B^k$ while $i_k:B^k\to B^{k-1}$ denotes the canonical projection. Either way gives rise to essentially the same topological realisation functor as can be seen by comparison of their values at representable presheaves. The recursive formula for $|\Theta_n[T]|_{\Theta_n}$ (as an iterated simplicial suspension of the one-point space, cf. the proof of \ref{realisation}) can also be found in \cite[Proposition 3]{Joy1}; the direct non-recursive formula is described in detail in \cite[Proposition 2.6]{mi}; both methods yield a $\CW$-ball of dimension equal to the number of edges of the representing level-tree $T$.\end{proof}

\begin{rmk}\label{mix}This topological realisation functor $|\!-\!|_{\Theta_n}:\widehat{\Theta}_n\to\Top$ mixes in a clever way simplicial and hemispherical combinatorics. Indeed, the presheaf represented by a linear level-tree of height $k$ realises to a $k$-ball with its hemispherical cell-decomposition, while the presheaf represented by a $1$-level-tree with $k$ edges realises to a euclidean $k$-simplex with its simplicial cell-decomposition. In general, the presheaf represented by a level-tree $T$ realises to a convex subset of a cube of dimension equal to the number of edges of $T$, cf. \cite[Proposition 2.6]{mi}. The left exactness of the realisation functor is related to a shuffle-formula which decomposes a product of representable presheaves $\Theta_n[S]\times\Theta_n[T]$ into a union of representable presheaves $\Theta_n[U]$, where $U$ runs through the set of level-trees that are obtained by shuffling $S$ and $T$ over the root. This formula generalises the decomposition of $\Delta[m]\times\Delta[n]$ into $\frac{(n+m)!}{n!m!}$ copies of $\Delta[n+m]$, cf. \cite[Proposition 2.8]{mi}. 

Proposition \ref{realisation} shows that the realisation functor for $\widehat{\Theta}_n$ factors through the realisation functor for $\widehat{\Delta^{\times n}}$ by means of an iterated diagonal $\delta_n:\Delta^{\times n}\to\Theta_n$, in other words, there is a natural isomorphism of functors $|\!-\!|_{\Theta_n}\cong|\!-\!|_{\Delta^{\times n}}\circ(\delta_n)^*$. This implies in particular that \emph{strict $n$-categories} may be realised topologically in two different, but homeomorphic ways. Indeed, the dense embedding $\Theta_n\inc\nCat$ defines a fully faithful nerve functor $\NN_{\Theta_n}:\nCat\to\widehat{\Theta}_n$; composing the latter with $|\!-\!|_{\Theta_n}:\widehat{\Theta}_n\to\Top$ yields a topological realisation for strict $n$-categories. Alternatively, a strict $n$-category may be realised by applying iteratively the simplicial nerve functor, and then realising the resulting $n$-simplicial nerve. The $n$-simplicial nerve functor is easily identified with the composite functor $(\delta_n)^*\circ\NN_{\Theta_n}$; therefore, the two realisations of a strict $n$-category coincide up to homeomorphism. However, the combinatorial structures of the two nerves are different; in particular, the $\Theta_n$-nerve is a full functor, while the $n$-simplicial nerve is \emph{not} a full functor for $n>1$. 

Despite this difference, strict $n$-categories can be characterised as special nerves using either the $n$-simplicial or the $\Theta_n$-nerve, and a definition of weak $n$-category may be deduced from both approaches. The $n$-simplicial approach leads to Simpson-Tamsamani's definition \cite{Tam} while the $\Theta_n$-approach leads to Joyal's definition \cite{Joy1} of a weak $n$-category. The iterated diagonal $\delta_n:\Delta^{\times n}\to\Theta_n$ induces a direct comparison of the two definitions, cf. \cite[Remarks 1.13/2.13]{mi}.\end{rmk}

\subsection{Pregeometric Reedy categories}\label{factor} An explicit non-recursive description of the operators in $\,\Theta_n$ or the dual operators in $\DD_n$ affords some combinatorial work, which has been done by several authors, cf. \cite{Ba1}, \cite{BS}, \cite{mi}, \cite{Joy1}, \cite{MZ}. We advertise the reader that a comparison of the different descriptions is a non-trivial task, some hints are given in Remark \ref{factor2} below. Later on, we explicitly use one further structure on $\,\Theta_n$, namely its \emph{Reedy structure}. We could derive this Reedy structure from \cite[Lemma 2.4]{mi} on the basis of Theorem \ref{dense}, which allows us to identify the iterated wreath product $\Theta_n$ with the $n$-th filtration of the category $\Theta$ studied in detail in \cite{mi}. It is however more in the spirit of the present text to rederive the Reedy structure on $\Theta_n$ by induction on $n$. Interestingly, a naive induction does not work. We have to replace the notion of a Reedy category by a slightly stronger notion, which is the starting point of a forthcoming joint work with Denis-Charles Cisinski, see \cite{BC}. Here, we discuss just that part of the structure that is needed to define inductively the Reedy structure on $\Theta_n$, cf. Remark \ref{geometric}. Let $\AA$ be an arbitrary small category.\vspace{1ex}

\begin{itemize}\item $\AA$ is called \emph{retractive} if any morphism of $\AA$ factors in a unique way as a retraction followed by a monomorphism, and if moreover for any object $A$ of $\AA$, the poset of retractions under $A$ is a \emph{lattice};

\item $\AA$ is called \emph{cellular} if for any object $A$ of $\AA$, the poset of subobjects of $A$ is the cell-poset of a finite regular $\CW$-complex (such posets are called \emph{CW-posets} by Bj\"orner \cite{Bjo}; they are in particular ranked posets);

\item $\AA$ is called \emph{pregeometric} if $\AA$ is at once retractive and cellular.\end{itemize}

Retractive categories $\AA$ have the following four properties: the identities are the only isomorphisms; all idempotents of $\AA$ split; any map of presheaves $\AA[A]\to X$ factors in a unique way as a retraction $\AA[A]\cto\AA[B]$ followed by a \emph{non-degenerate} map $\AA[B]\to X$ (i.e. one that does not factor through a non-identity retraction); any finite product of representable presheaves is the union of its representable subobjects. (The last two properties follow from the lattice property of retractions, while the first two properties follow from the factorisation property of $\AA$). In a cellular category $\AA$, any object $A$ has a well-defined dimension $\dim_\AA(A)$, namely the rank of the poset of subobjects of $A$. Non-identity monomorphisms of $\AA$ \emph{increase} this dimension. Recall \cite{Hir},\cite{Ho} that a \emph{Reedy category} is a quadruple $(\AA,\AA_+,\AA_-,\dim_\AA)$ where $\AA$ is a small category whose objects are graded by the dimension-function $\dim_\AA$. Moreover, $\AA$ contains subcategories $\AA_+$ resp. $\AA_-$ whose non-identity morphisms increase resp. decrease dimension, and any morphism of $\AA$ factors in a unique way as a morphism in $\AA_-$ followed by a morphism in $\AA_+$.

It follows that any pregeometric category $\AA$ has a canonical Reedy structure $(\AA,\AA_+,\AA_-,\dim_\AA)$ where $\AA_+$ (resp. $\AA_-$) contains all monomorphisms (resp. retractions) of $\AA$. A pregeometric category $\AA$ equipped with this Reedy structure will be called a \emph{pregeometric Reedy category}. The simplex category $\Delta$ with its canonical Reedy structure is a pregeometric Reedy category; the only point needing a comment is the lattice property of the retractions; the latter expresses the fact that the poset of ordered partitions of $[n]$ is a lattice for any $n$. Usually, the monomorphisms (resp. retractions) of $\Delta$ are called \emph{face} (resp. \emph{degeneracy}) operators; we shall adopt the same terminology for an arbitrary pregeometric Reedy category $\AA$.

Any pregeometric Reedy category $\AA$ has a canonical embedding in the category of finite regular $\CW$-complexes and cellular maps; indeed, for any object $A$ of $\AA$, the simplicial nerve of the subobject-poset $\AA_+/A$ represents the \emph{barycentric subdivision} of the regular $\CW$-complex $C_A$ whose cell-poset is $\AA_+/A$. Therefore, the realisation $|\NN_\Delta(\AA_+/A)|_\Delta$ is canonically homeomorphic to $C_A$. Each operator $A\to B$ in $\AA$ induces thus (by factoring it into a degeneracy followed by a face operator) a functor $\AA_+/A\to\AA_+/B$ and hence a cellular map $C_A\to C_B$. Left Kan extension of this embedding along the Yoneda-embedding $\AA\inc\widehat{\AA}$ defines a colimit preserving functor $|\!-\!|_\AA:\widehat{\AA}\to\Top$ with values in the category of $\CW$-complexes. Indeed, for any presheaf $X$ on $\AA$, the topological space $|X|_\AA$ decomposes as a (set-theoretically) disjoint union of ``cell-interiors'', indexed by the non-degenerate elements of $X$, exactly as in the case of a simplicial set. 

A pregeometric Reedy category $\AA$ will be called \emph{(locally) flat} if $|\!-\!|_\AA:\widehat{\AA}\to\Top$ preserves finite (connected) limits. In particular, $\AA$ is flat if and only if $\AA$ is locally flat and the terminal presheaf $\star_\AA$ realises to the one-point space. The latter happens for instance if $\AA$ contains a terminal object without proper subobjects. The next proposition establishes property (3) of the introduction to this section.

\begin{prp}\label{pregeometric}For any flat pregeometric Reedy category $\AA$, the wreath product $\Delta\wr\AA$ is again a flat pregeometric Reedy category with dimension-function$$\dim_{\Delta\wr\AA}((A_1,\dots,A_m))=m+\dim_\AA(A_1)+\cdots+\dim_\AA(A_m).$$In particular, the iterated wreath product $\,\Theta_n$ is a flat pregeometric Reedy category with dimension-function $\dim_{\Theta_n}(T)$ given by the number of edges of the level-tree $T$.\end{prp}

\begin{proof}We have to show that $\Delta\wr\AA$ is retractive, cellular and flat, with the indicated dimension-function. The second statement follows then by induction on $n$.

Let $(\phi;\phi_1,\dots,\phi_m):(A_1,\dots,A_m)\to(B_1,\dots,B_n)$ be an operator in $\Delta\wr\AA$. Such an operator is a retraction precisely when $\phi$ is a retraction, and each $\phi_k$ is either the trivial morphism $\AA[A_k]\to\star_\AA$ or a retraction $\AA[A_k]\to\AA[B_{k'}]$ in $\widehat{\AA}$. For a general operator $(\phi;\phi_1,\dots,\phi_m)$, factor $\phi$ as a retraction $\psi$ followed by a monomorphism $\rho$. Then $(\phi;\phi_1,\dots,\phi_m)$ factors as $(\psi;\psi_1,\dots,\psi_m)$ followed by $(\rho;\rho_1,\dots,\rho_{m'})$ where each $\psi_k$ is either trivial or an identity. Since $\AA$ is retractive, each $\rho_k:\AA[A_k]\to\AA[B_{\rho(k-1)+1}]\times\cdots\times\AA[B_{\rho(k)}]$ factors as a retraction $\sg_k:\AA[A_k]\to\AA[A'_k]$ followed by a non-degenerate map $\tau_k:\AA[A'_k]\to\AA[B_{\tau(k-1)+1}]\times\cdots\times\AA[B_{\tau(k)}]$ (i.e. $\tau_k$ does not factor through any non-identity retraction). The retractions $\sg_k$ assemble into a retraction $(id_{[m']};\sg_1,\dots\sg_{m'})$ in $\Delta\wr\AA$. Since $\AA$ is retractive, the non-degenerate map $\tau_k$ is monic in $\widehat{\AA}$ because otherwise the target of $\tau_k$ would not be the union of its representable subobjects. Therefore, we get the required factorisation of $(\phi;\phi_1,\dots,\phi_m)$ into a retraction $(id_{[m']};\sg_1,\dots\sg_{m'})(\psi;\psi_1,\dots,\psi_m)$ followed by a monomorphism $(\rho;\tau_1,\dots,\tau_{m'})$. Assume there is another such factorisation. Then, it follows from the above that the two retraction-parts are the same, and hence the two monic parts are the same too. The lattice property of the retractions in $\Delta\wr\AA$ follows from the lattice property of the retractions in $\Delta$ and the lattice property of the retractions in $\AA$.

For the cellularity of $\Delta\wr\AA$, we use diagram (\ref{suspension2}) of the proof of Proposition \ref{realisation}, which identifies $(\Delta\wr\AA)[(A_1,\dots,A_m)]$ with a certain subobject of the cartesian product $(\Delta\wr\AA)[(A_1)]\times\cdots\times(\Delta\wr\AA)[(A_m)]$. We have to show that the realisation functor $|\!-\!|_{\Delta\wr\AA}$ of Proposition \ref{realisation} endows $|(\Delta\wr\AA)[(A_1,\dots,A_m)]|_{\Delta\wr\AA}$ with the structure of a regular $\CW$-complex with cell-poset $(\Delta\wr\AA)_+/(A_1,\dots,A_m)$. This is clear for $m=1$, since for any object $A$, the unreduced suspension of $|\AA[A]|_\AA$ may be obtained by applying the unreduced suspension cell-wise and gluing the suspended cells together like in the original cell; the resulting cell-poset of the unreduced suspension of $|\AA[A]|_\AA$ is then precisely $(\Delta\wr\AA)_+/(A)$. Therefore, we get a canonical product $\CW$-structure on $|(\Delta\wr\AA)[(A_1)]\times\cdots\times(\Delta\wr\AA)[(A_m)]|_{\Delta\wr\AA}$. The right vertical map of diagram (\ref{suspension2}) maps this product cellularly to the $m$-cube $(\Delta^1)^{\times m}$ endowed with its canonical product $\CW$-structure. The latter admits a regular $\CW$-subdivision by $m!$ copies of the standard $m$-simplex $\Delta^m$; this subdivision lifts in a unique way to a regular $\CW$-subdivision of the product $|(\Delta\wr\AA)[(A_1)]\times\cdots\times(\Delta\wr\AA)[(A_m)]|_{\Delta\wr\AA}$. It is then straightforward to check that the cells of this regular $\CW$-structure correspond canonically to the representable subobjects of $(\Delta\wr\AA)[(A_1)]\times\cdots\times(\Delta\wr\AA)[(A_m)]$. In particular, the representable subobject $(\Delta\wr\AA)[(A_1,\dots,A_m)]$ realises to a regular $\CW$-complex with cell-poset $(\Delta\wr\AA)_+/(A_1,\dots,A_m)$.

The flatness of $\Delta\wr\AA$ follows from Proposition \ref{realisation}, since the realisation functor defined there coincides up to homeomorphism with the realisation functor derived from the cellular structure of $\Delta\wr\AA$. The formula for the dimension-function expresses the fact that $(\Delta\wr\AA)[(A_1,\dots,A_m)]$ and $(\Delta\wr\AA)[(A_1)]\times\cdots\times(\Delta\wr\AA)[(A_m)]$ realise to $\CW$-balls of the same dimension.\end{proof}

\begin{rmk}\label{factor2}The degeneracy operators of $\Theta_n$ admit a nice description in terms of level-trees. Indeed, simplicial degeneracy operators $[k]\to[l]$ correspond bijectively to ``subtree''-inclusions of a corolla with $l$ leaves into a corolla with $k$ leaves; this implies by induction on $n$ that degeneracy operators $S\to T$ in $\Theta_n$ correspond bijectively to ``subtree''-inclusions of $T$ into $S$ on the basis of the following definition, cf.  Section \ref{level-tree} and \cite[Definition 2.3]{mi}: a \emph{subtree} $T$ of a level-tree $S$ is a subgraph of $S$ sharing the same root, such that for each vertex $v$ of $T$ the edges of $T$ incoming in $v$ form a \emph{connected} subset of the set of edges of $S$ incoming in $v$.

The face operators of $\Theta_n$ may be divided into \emph{inner} and \emph{outer} face operators in such a way that each face operator factors uniquely as an inner followed by an outer face operator. This is again proved by induction on $n$. A simplicial face operator is called inner (resp. outer) if it preserves minimal and maximal elements (resp. if $\phi(i+1)=\phi(i)+1$ for all $i$). Then by induction on $n$, a face operator $(\phi;\phi_1,\dots,\phi_m)$ of $\Theta_n$ is called inner (resp. outer) if $\phi$ is an inner face operator of $\Delta$ and each projection of each $\phi_k:\Theta_{n-1}[S_k]\to\Theta_{n-1}[T_{\phi(k-1)+1}]\times\cdots\times\Theta_{n-1}[T_{\phi(k)}]$ is an inner face operator of $\Theta_{n-1}$ (resp. if $\phi$ is an outer face operator of $\Delta$ and each $\phi_k$ is either trivial or an outer face operator of $\Theta_{n-1}$).

Inner face operators also admit a nice description in terms of level-trees. Indeed, an inner face operator $\phi:S\to T$ determines, and is determined by, a family $(T_v)_{v\in in(S)}$ of subtrees of $T$ indexed by the input vertices $v\in in(S)$ such that (i) the height of $v$ equals the height of $T_v$; (ii) the union of the $T_v$ is $T$; and (iii) for any pair of input vertices $(v_1,v_2)\in in(S)\times in(S)$ for which the edge-paths to the root meet at height $k$, the $k$-level-truncations of $T_{v_1}$ and $T_{v_2}$ coincide with the intersection $T_{v_1}\cap T_{v_2}$. The composite of two such operators \begin{diagram}S&\rTo^{(T_v)_{v\in in(S)}}&T&\rTo^{(U_w)_{w\in in(T)}}U\end{diagram}is obtained by assigning to $v\in in(S)$ the union of all $U_w$ for $w\in in(T_v)$.

The subcategory of $\Theta_n$ containing the outer face operators can be identified with the image of the free functor from $n$-graphs to $n$-categories; in other words, outer face operators $S\to T$ correspond to maps of $n$-graphs $S_*\to T_*$, cf. Section \ref{level-tree}. In particular, the subcategory of outer face operators comes equipped with a subcanonical Grothendieck topology (where the covering sieves are epimorphic families); the latter can be used to \emph{characterise} nerves of $n$-categories as those presheaves on $\Theta_n$ which become sheaves when restricted to the subcategory of outer face operators, cf. \cite[Remark 1.13]{mi}. For $n=1$, this recovers the well known Grothendieck-Segal characterisation of nerves of categories. This also explains the role of outer face operators in the definition of a Segal-fibration, cf. Definition \ref{Segalfibration}.

In \cite{mi}, face operators are called \emph{level-preserving operators}, outer face operators are called \emph{immersions}, and operators of the subcategory generated by inner face operators and degeneracy operators are called \emph{covers}. Joyal \cite{Joy1} calls the dual of such a cover an \emph{open disk map}. Batanin \cite{Ba1} calls inner face operators \emph{tree-diagrams}. 

Batanin and Street \cite{BS} define a category $\Omega_n$ which is closely related to Joyal's category of finite combinatorial $n$-disks and open disk maps; however, $\Omega_n$ contains in addition for each level-tree $S$ of height $< n$ \emph{degenerate level-trees} $Z^kS$ for $0<k\leq n-ht(S)$, to the effect that those open disk maps which are dual to degeneracy operators get a slightly different description in $\Omega_n$. Each of the categories $\Omega_n$ is characterised by a universal property; indeed, $\Omega_n$ is the precise $n$-categorical analog of $\Omega_1\cong\Delta_{alg}=\Delta\cup\{[-1]\}$, where $[-1]$ denotes the empty ordinal.

The assembly functor $\gamma_n:\Theta_n\to\Gamma$ coincides after restriction to covers and dualisation (up to the different ways of representing degeneracies) with Batanin's functor $\Omega_n\to\Omega^s\subset\Gamma^\op$, cf. \cite{Ba2}. Indeed, the assembly functor associates to each operator $S\to T$ its ``trace'' on the vertices of height $n$; in particular, the value $\gamma_n(T)$ is $\kk$ iff $T$ has $k$ vertices of height $n$. It follows then from the combinatorial description of a cover $\phi:S\to T$ given above that the induced operator $\gamma_n(\phi):\gamma_n(S)\to\gamma_n(T)$ is a $\Gamma$-operator of a special type: it ``covers'' its target, or equivalently, the dual operator is a map of finite pointed sets which restricts away from the base-point, i.e. which belongs to Batanin's subcategory $\Omega^s$ of $\Gamma^\op$.\end{rmk}

\begin{prp}\label{classifying}For each $n\geq 1$, the assembly functor $\gamma_n:\Theta_n\to\Gamma$ induces a geometric morphism of presheaf toposes $(\gamma_n)_*:\widehat{\Theta}_n\rlto\widehat{\Gamma}:(\gamma_n)^*$ which classifies the quotient of the generic combinatorial $n$-disk by its boundary.\end{prp}

\begin{proof}For $n=1$, the generic combinatorial $1$-disk is the standard $1$-simplex $\Delta[1]$ of the classifying topos for combinatorial $1$-disks. Since the $\Gamma$-set $\Gamma[\underline{1}]$ is the generic pointed object of the classifying topos for pointed objects, we have to show that $\gamma_1^*\Gamma[\underline{1}]=\Delta[1]/\partial\Delta[1]$. This follows from the fact that distinct simplicial operators $[m]\to[1]$ are identified under $\gamma_1:\Delta\to\Gamma$ if and only if they factor through $[0]$.

By induction on $n$, the generic combinatorial $n$-disk of the classifying topos for combinatorial $n$-disks is represented by the linear level-tree $\bar{n}$ of height $n$, cf. the proof of Theorem \ref{disk}. We have to show that $\gamma_n^*\Gamma[\underline{1}]=\Theta_n[\bar{n}]/\partial\Theta_n[\bar{n}]$ where the boundary $\partial\Theta_n[\bar{n}]$ is derived from the Reedy structure on $\Theta_n$, cf. Proposition \ref{pregeometric}.

By induction on $n$, the value $\gamma_n(T)$ at a level-tree $T$ with $k$ vertices of height $n$ is the $k$-element set $\kk$. In particular, $\gamma_n(T)$ is the null-object of $\Gamma$ if and only if $ht(T)<n$. Moreover, $\gamma_n(T)=\underline{1}$ if and only if $\bar{n}$ embeds as a subtree in $T$ in a unique way, i.e. (cf. Remark \ref{factor2}) iff there is a uniquely determined degeneracy operator $\rho:T\to\nn$; in particular $\rho\in\Theta_n[\nn]$. Since $\partial\Theta_n[\bar{n}]$ consists of all operators $T\to\nn$ which factor through a level-tree with $<n$ edges, it remains to be shown that distinct operators $\rho_1,\rho_2:T\dto\bar{n}$ have same image under $\gamma_n$ if and only if they factor through a level-tree with $<n$ edges. 

If $\rho_1$ and $\rho_2$ factor through a level-tree with $< n$ edges then $\gamma_n(\rho_1)$ and $\gamma_n(\rho_2)$ factor through the null-object of $\Gamma$ and therefore coincide. Conversely, if $\gamma_n(\rho_1)=\gamma_n(\rho_2)$, then $(\gamma_1\wr\gamma_{n-1})(\rho_1)=(\gamma_1\wr\gamma_{n-1})(\rho_2)$ since the assembly functor $\alpha:\Gamma\wr\Gamma\to\Gamma$ of Proposition \ref{semi-additive} is faithful. Now, assume that $T=([m];S^1,\dots,S^m)$, and that $\rho_1=(\phi;\phi_1,\dots,\phi_m)$ and $\rho_2=(\psi;\psi_1,\dots,\psi_m)$, in particular $\gamma_1(\phi)=\gamma_1(\psi)$. If $\phi\not=\psi$, this implies that $\phi$ and $\psi$ factor through $[0]$ and we are done; if $\phi=\psi:[m]\to[1]$ is a degeneracy operator, then all $\phi_j$ and $\psi_j$ for $j\not=i$ are trivial, where $\phi(i-1)=0$ and $\phi(i)=1$, while $\phi_i\not=\psi_i$ and $\gamma_{n-1}(\phi_i)=\gamma_{n-1}(\psi_i)$. By induction hypothesis, this implies that $\phi_i$ and $\psi_i$ factor through a level-tree with $<(n-1)$ edges; therefore, $\rho_1$ and $\rho_2$ factor through a level-tree with $< n$ edges.\end{proof}

\begin{cor}\label{homeo}For any $\Gamma$-space $A$ and any $n\geq 1$, the geometric realisation $|\gamma_n^*A|_{\Theta_n}$ is homeomorphic to the $n$-th space $\underline{A}(S^n)$ of the Segal spectrum of $A$.\end{cor}

\begin{proof}The left Kan extension of the flat functor $\Theta^-_n:\Theta_n\inc\widehat{\Theta}_n\to\Top$ along the assembly functor $\gamma_n:\Theta_n\to\Gamma$ yields a flat and hence left exact functor $\gamma_{n!}\Theta^-_n:\Gamma\to\Top$; the latter is determined by its value at $\underline{1}$ which by adjunction is $$|\gamma_n^*\Gamma[\underline{1}]|_{\Theta_n}=|\Theta_n[\bar{n}]/\partial\Theta_n[\bar{n}]|_{\Theta_n}=B^n/\partial B^n=S^n.$$Therefore, $\gamma_{n!}\Theta^-_n:\Gamma\to\Top$ is the functor which takes $\kk$ to $(S^n)^k$. This implies for any $\Gamma$-space $A$ the identifications\begin{align*}\label{gamma2}|\gamma_n^*A|_{\Theta_n}=(\gamma_n^*A)\otimes_{\Theta_n}\Theta^-_n=A\otimes_\Gamma\gamma_{n!}\Theta_n^-=A\otimes_\Gamma (S^n)^-=\underline{A}(S^n).\end{align*}\end{proof}

\begin{rmk}The preceding proposition establishes property (4) of the introduction to this section; its corollary relates the combinatorics of $\Theta_n$-sets in an interesting way to stable homotopy theory. We are grateful to Bj\o rn Dundas for pointing out to us the following more concrete way of understanding this relationship.

  The dual $\Gamma^{op}$ of $\Gamma$ is the category of finite pointed sets and embeds therefore canonically in $\Sets$. The classifying property of Segal's functor $\gamma:\Delta\to\Gamma$ is then just a fancy way of saying that its dual $\gamma^\op:\Delta^\op\to\Gamma^\op\subset\Sets$ is precisely the simplicial circle $\Delta[1]/\partial\Delta[1]$. Similarly, the classifying property of the assembly functor $\gamma_n:\Theta_n\to\Gamma$ is just a fancy way of saying that its dual $\gamma_n^\op:\Theta_n^\op\to\Gamma^\op\subset\Sets$ is precisely the canonical $\Theta_n$-set-model $\Theta_n[\bar{n}]/\partial\Theta_n[\bar{n}]$ of the $n$-sphere. This last statement can be related to some classical constructions with $\Gamma$-spaces. Indeed, consider the following commutative diagram of functors\begin{diagram}[noPS,small,silent]\Delta\times\cdots\times\Delta&\rTo^{\delta_n}&\Delta\wr\cdots\wr\Delta&&\\\dTo&&\dTo&\rdTo^{\gamma_n}&\\\Gamma\times\cdots\times\Gamma&\rTo^{\delta^\Gamma_n}&\Gamma\wr\cdots\wr\Gamma&\rTo^\alpha&\Gamma\end{diagram}where the vertical functors are induced by $\gamma:\Delta\to\Gamma$ and $\delta_n^\Gamma$ is defined in complete analogy to $\delta_n$, cf. Definition \ref{diagonal}. In particular, the composite of $\delta^\Gamma_n$ with the iterated assembly functor $\alpha$ yields the $n$-fold product map $m_n:\Gamma^{\times n}\to\Gamma$. Proposition \ref{realisation} shows that the $\Theta_n$-set $\gamma_n^\op:\Theta_n^\op\to\Sets$ can be realised by pulling back along $\delta_n$ and realising the resulting $n$-simplicial set. The commutativity of the diagram above shows then that this $n$-simplicial set is precisely the \emph{$n$-fold smash product} of the simplicial circle, and hence realises to an $n$-sphere.

Lydakis \cite{Lyd} shows that the smash-product of $\Gamma$-spaces, defined by left Kan extension along $m_2:\Gamma\times\Gamma\to\Gamma$, has homotopical meaning insofar as monoids with respect to this smash-product (so-called $\Gamma$-rings) are models for connective ring spectra, cf. Remark \ref{connectivespectra}. An essential ingredient in his proof is the homotopical analysis of the so-called \emph{assembly map} $A\wedge B\to A\circ B$ for $\Gamma$-spaces $A$ and $B$. This assembly map is closely related to the diagonal $\delta_2^\Gamma:\Gamma\times\Gamma\to\Gamma\wr\Gamma$, and therefore, the fact that the assembly map is a weak equivalence (if $A$ or $B$ is cofibrant), is closely related to Proposition \ref{realisation}.\end{rmk}

\subsection{Embedding functors}\label{embedding}The category $\Theta_n$ embeds in $\Theta_{n+1}$ as the full subcategory spanned by the level-trees of height $\leq n$. This embedding is compatible with the embedding of Theorem \ref{dense} in the sense that we have commutative squares\begin{diagram}[noPS,small,silent]\Theta_n&\rTo&\nCat\\\dTo&&\dTo\\\Theta_{n+1}&\rTo&\nplusCat.\end{diagram}The embedding $i_n:\Theta_{n}\to\Theta_{n+1}$ is also compatible with the realisation functors of Corollary \ref{flat} in the sense that we have commutative triangles\begin{diagram}[noPS,small,silent]\widehat{\Theta}_n&\rTo^{|\!-\!|_{\Theta_n}}&\Top\\\dTo^{i_{n!}}&\ruTo_{|\!-\!|_{\Theta_{n+1}}}&\\\widehat{\Theta}_{n+1}&&\end{diagram}where as always $i_{n!}$ denotes left Kan extension along $i_n$.

\subsection{Suspension functors}\label{suspensionfunctor}The inductive definition $\Theta_{n+1}=\Delta\wr\Theta_n$ defines a \emph{second} completely different functor $$\sg_n:\Theta_n\to\Theta_{n+1},$$ given on objects by $S\mapsto(S)$ and on operators by $\rho\mapsto(id_{[1]};\rho)$. The level-tree $(S)$ is the level-tree $S$ shifted upwards one level by adjunction of an extra root-edge. We call $\sg_n$ a \emph{suspension functor}; this functor has already been considered by Joyal \cite{Joy1} and also (in dual form, cf. Remark \ref{factor2}) by Batanin \cite{Ba2}. It follows from the definitions that we have commutative triangles for $n\geq 0$\begin{diagram}[noPS,small,silent]\Theta_n&\rTo^{\gamma_n}&\Gamma\\\dTo^{\sg_n}&\ruTo_{\gamma_{n+1}}&&\\\Theta_{n+1}&&&\end{diagram}where $\Theta_0$ denotes the terminal category and $\gamma_0:\Theta_0\to\Gamma$ takes the unique object of $\Theta_0$ to $\underline{1}$. Let us mention that $\Theta_0$ has a completely different meaning in \cite{mi}, where it stands for the subcategory of outer face operators of the union $\Theta$ of the $\Theta_n$.

The suspension functors $\sg_n$ are compatible with the \emph{unreduced} topological suspension functor $S$ as expressed by the following (pseudo)commutative squares\begin{diagram}[noPS,small,silent]\Theta_n&\rTo^{\Theta^-_n}&\Top\\\dTo^{\sg_n}&&\dTo_S\\\Theta_{n+1}&\rTo^{\Theta^-_{n+1}}&\Top\end{diagram}where the flat functors $\Theta_n^-:\Theta_n\inc\widehat{\Theta}_n\overset{|\!-\!|_{\Theta_n}}{\lra}\Top$ can indifferently be defined using either Proposition \ref{realisation} or Corollary \ref{flat} or Proposition \ref{pregeometric}.

\section{Reduced $\Theta_n$-spaces and $n$-fold loop spaces}

\begin{dfn}A $\Theta_n$-space $X:\Theta_n^\op\to\Top$ is \emph{reduced} if the value of $X$ at each level-tree of height $<n$ is the one-point space. A \emph{$\Theta$-spectrum} $(X_n,s_n)_{n\geq 0}$ is a sequence of reduced $\Theta_n$-spaces $X_n$ and natural transformations $s_n:X_n\to \sg_n^*X_{n+1}$.\end{dfn}

\begin{prp}\label{thetaspectra}The Segal spectrum functor from $\Gamma$-spaces to topological spectra factors through the category of $\Theta$-spectra via $A\mapsto (\gamma_n^*A)_{n\geq 0}$.\end{prp}
\begin{proof}For any $\Gamma$-space $A$, the canonical isomorphisms $s_n:\gamma_n^*A\to\sigma_n^*\gamma_{n+1}^*A$ define a $\Theta$-spectrum $(\gamma_n^*A,s_n)_{n\geq 0}$. Any $\Theta$-spectrum $(X_n,s_n)_{n\geq 0}$ defines a topological spectrum $(|X_n|_{\Theta_n})_{n\geq 0}$ with the following structural maps: according to the proof of Proposition \ref{realisation}, the geometric realisation $|X_{n+1}|_{\Theta_{n+1}}$ may be identified with the realisation of the simplicial space $r_{\Theta_n}X_{n+1}=|\delta_{\Theta_n}^*X_{n+1}|_{\Theta_n}$; this simplicial space is \emph{reduced}, and its \emph{underlying space} is precisely $|\sg_n^*X_{n+1}|_{\Theta_n}$, cf. Section \ref{redsimpspaces}; we get thus a canonical map $|\sg_n^*X_{n+1}|_{\Theta_n}\to\Omega|X_{n+1}|_{\Theta_{n+1}}$ and hence, by precomposition with $|s_n|_{\Theta_n}$, the structural maps $|X_n|_{\Theta_n}\to \Omega|X_{n+1}|_{\Theta_{n+1}}$ of a topological spectrum. It follows from Corollary \ref{homeo} that the composite of the two functors above is isomorphic to the Segal spectrum functor $\Phi:\Top_{red}^{\Gamma^\op}\to\Spt$ of Section \ref{Gammaspaces}.\end{proof}

Each abelian group $\pi$ gives rise to a discrete $\Gamma$-space $H\pi$ (also called a $\Gamma$-set) defined by $(H\pi)(\kk)=\pi^k$ where the $\Gamma$-operations are induced by the abelian structure of $\pi$. This $\Gamma$-space $H\pi$ is cofibrant-stably fibrant for the injective model structure of Bousfield and Friedlander; its associated Segal spectrum $\Phi(H\pi)$ is an explicit model for the \emph{Eilenberg-MacLane spectrum} of the abelian group $\pi$, cf. \cite[1.5]{Sch}. In particular, the $n$-th space $\underline{H\pi}(S^n)$ has the homotopy type of an \emph{Eilenberg-MacLane space} of type $K(\pi,n)$. The preceding proposition implies thus 

\begin{cor}\label{EM}For each abelian group $\pi$ and each $n$, the reduced $\Theta_n$-set $\gamma^*_n(H\pi)$ realises to an Eilenberg-MacLane space of type $K(\pi,n)$.\end{cor}

\begin{rmk}\label{geometric}The methods of \cite[chapter  3]{mi} yield a Quillen model structure on $\widehat{\Theta}_n$ together with a monoidal Quillen equivalence $|\!-\!|_{\Theta_n}:\widehat{\Theta}_n\lrto\Top:\Sing_{\Theta_n}.$ In joint work with Cisinski \cite{BC}, a new proof of this result is obtained along with a far reaching generalisation. Indeed, a Reedy category $\AA$ is called \emph{geometric} in loc. cit. if it is a pregeometric Reedy category in the sense of Section \ref{factor} with the further property that $\widehat{\AA}$ admits \emph{functorial cylinder-objects} in Quillen's sense \cite{Q}, where the cofibrations (resp. weak equivalences) are the injective maps (resp. those maps of presheaves that realise to weak homotopy equivalences). One of the key results of \cite{BC} says that for any (locally flat) geometric Reedy category $\AA$, this defines a proper Quillen model structure on $\widehat{\AA}$, cofibrantly generated by sphere- and horn-inclusions, together with a (monoidal) Quillen equivalence $|\!-\!|_\AA:\widehat{\AA}\lrto\Top/|\!\star_{\widehat{\AA}}\!|_\AA:\Sing_\AA$. 

It follows from Proposition \ref{pregeometric} that the iterated wreath-product $\Theta_n$ is a flat geometric Reedy category; this means thus that the presheaf topos $\widehat{\Theta}_n$ has a full-fledged homotopy theory, including a well-behaved combinatorial concept of fibration, and a homotopy-invariant notion of Eilenberg-MacLane object of type $K(\pi,n)$. The reduced $\Theta_n$-sets $\gamma_n^*(H\pi)$ are particularly nice such objects. Surprisingly, the resulting $\CW$-complexes $|\gamma_n(H\pi)|_{\Theta_n}$ are comparatively small and permit a new, combinatorially interesting, ``cochain calculus''.\end{rmk}

Before we study in more detail these models for Eilenberg-MacLane spaces, we come to the main result of this article and endow the category $\Top_{red}^{\Theta_n^\op}$ of reduced $\Theta_n$-spaces with a Quillen model structure in such a way that it becomes a model for the derived image of the $n$-fold loop functor $\Omega^n:\Top_*\to\Top_*$. We will use Lemma \ref{nfold} and define first a commutative triangle of Quillen adjunctions\begin{gather}\begin{diagram}[noPS,p=1mm]\Top_*&\pile{\lTo^{|\!-\!|_{\Theta_n}}\\\rTo_{\Omega^n_{Seg}}}&\Top_{red}^{\Theta_n^\op}\\&\pile{\rdTo^{\Omega^n}\\\luTo_{\Sigma^n}}&\dTo^U\uTo_L\\&&\Top_*.\end{diagram}\end{gather}This triangle generalises triangle (\ref{1folddiag}) of Section \ref{redsimpspaces}. In particular, for any based space $(X,*)$, the $n$-fold Segal loop space $\Omega^n_{Seg}(X)$ is defined by$$\Omega^n_{Seg}(X)(T)=\Top(|\Theta_n[T]|_{\Theta_n},|i_{n-1!}i_{n-1}^*\Theta_n[T]|_{\Theta_n}),(X,*))$$with the usual (compact-open) topology, where $i_{n-1}:\Theta_{n-1}\to\Theta_n$ denotes the canonical embedding of Section \ref{embedding}. The underlying-space functor $U$ is given by $U(Y)=Y(\bar{n})$, where $\bar{n}$ denotes the linear level-tree of height $n$. Since $|\Theta_n[\bar{n}]|_{\Theta_n}=B^n$ and $|i_{n-1!}i_{n-1}^*\Theta_n[\bar{n}]|_{\Theta_n}=|\partial\Theta_n[\bar{n}]|_{\Theta_n}=\partial B^n$, cf. Remark \ref{mix}, the composite functor $U\Omega_{Seg}^n$ may be identified with the classical $n$-fold loop functor $\Omega^n$.

Since $\Theta_n$ is a Reedy category by Proposition \ref{pregeometric}, there is a canonical Reedy model structure on $\Top_*^{\Theta_n^\op}$ with pointwise weak equivalences, cf. \cite{Hir}, \cite{Ho}. The latter restricts to the bireflective subcategory $\Top_{red}^{\Theta_n^\op}$ of reduced $\Theta_n$-spaces by Remark \ref{specialtransfer} in such a way that the above-defined adjoint pair $(|\!-\!|_{\Theta_n},\Omega^n_{Seg})$ is a Quillen adjunction. A map of reduced $\Theta_n$-spaces is a \emph{realisation weak equivalence} if the left derived functor of $|\!-\!|_{\Theta_n}$ takes it to an isomorphism in $\Ho(\Top_*)$.

\begin{thm}\label{main}The Reedy model category of reduced $\,\Theta_n$-spaces admits a left Bous-field localisation with respect to realisation weak equivalences. The localised model category is a model for $n$-fold loop spaces. In particular, for any cofibrant-fibrant reduced $\,\Theta_n$-space $X$, the canonical map $U(X)\to\Omega^n|X|_{\Theta_n}$ is a weak equivalence.\end{thm}

\begin{proof}We argue by induction on $n$. The case $n=1$ is Theorem \ref{Segalfibrant}. The inductive step breaks into three parts: the existence of the asserted localisation, the property that the underlying-space functor is homotopy-right-conservative, and the property that the canonical map $U(X)\to\Omega^n|X|_{\Theta_n}$ is a weak equivalence at cofibrant-fibrant objects (which is equivalent to $\Omega^n_{Seg}$ being a homotopy-coreflection).

Recall from Proposition \ref{realisation} that the realisation functor $|\!-\!|_{\Theta_n}$ factors through the functor $\delta_{\Theta_{n-1}}^*:\Top_{red}^{\Theta_n^\op}\to(\Top_{red}^{\Theta^\op_{n-1}})_{red}^{\Dop}$; the suspension functor $\sg_{n-1}$ induces a functor $\sg_{n-1}^*:\Top_{red}^{\Theta_n^\op}\to\Top_{red}^{\Theta_{n-1}^\op}$ which commutes with the underlying-space functors. We already observed (in the proof of Proposition \ref{thetaspectra}) that the following triangle commutes up to isomorphism\begin{gather}\label{rs}\begin{diagram}[small,noPS,silent]\Top_{red}^{\Theta_n^\op}&\rTo^{\delta^*_{\Theta_{n-1}}}&(\Top_{red}^{\Theta_{n-1}^\op})_{red}^{\Dop}\\\dTo^{\sg_{n-1}^*}&\ldTo_U&\\\Top_{red}^{\Theta_{n-1}^\op}&&\end{diagram}\end{gather}

Inductively, a \emph{Segal-fibration} of reduced $\Theta_n$-spaces $f:X\to Y$ is defined to be a Reedy-fibration such that the following two conditions hold:\begin{itemize}\item[(a)] denote by $G_\FF$ the subobject of $\Theta_n[T]$ generated by a family $\FF$ of outer face operators $T_1\to T,\dots,T_m\to T$, which is epimorphic in the subcategory of outer face operators, cf. Remark \ref{factor2}; then for any such $G_\FF$, the induced map $X(T)\to Y(T)\times_{Y(G_\FF)}X(G_\FF)$ has to be a weak equivalence;\item[(b)]$\delta^*_{\Theta_{n-1}}(f)$ resp. $\sg_{n-1}^*(f)$ fulfill condition (b) for a Segal-fibration of reduced simplicial objects in $\Top_{red}^{\Theta_{n-1}^\op}$ resp. of reduced $\Theta_{n-1}$-spaces.\end{itemize}Now we apply Proposition \ref{localisation}. The existence of a suitable set of realisation-trivial cofibrations characterising Segal-fibrations among Reedy-fibrations follows by induction on $n$, since $\delta^*_{\Theta_{n-1}}$ and $\sg_{n-1}^*$ have left adjoints. It remains to be shown that any realisation-trivial Segal-fibration $f:X\to Y$ is a trivial Reedy-fibration. For this, consider a level-tree $T$ with root-valence $m$; then, $T$ is a bouquet of $m$ level-trees $T_1,\dots,T_m$ glued together at the root. The isomorphism of $n$-graphs $(T_1)_*\vee\cdots\vee(T_m)_*\cong T_*$ defines an epimorphic family of outer face operators $T_k\to T,\,1\leq k\leq m$; condition (a) of a Segal-fibration induces then the following homotopy cartesian square:\begin{diagram}[small,noPS]X(T)&\rTo^{f(T)}&Y(T)\\\dTo&&\dTo\\X(T_1)\times\cdots\times X(T_m)&\rTo^{f(T_1)\times\cdots\times f(T_m)}&Y(T_1)\times\cdots\times Y(T_m)\end{diagram}It follows from the inductive definition of a Segal-fibration that the induced maps $\sg^*_{n-1}(f)$ and $\delta^*_{\Theta_{n-1}}(f)$ are Segal-fibrations. Since $f$ is a realisation weak equivalence, $\delta^*_{\Theta_{n-1}}(f)$ is a realisation weak equivalence as well, and the analog of Theorem \ref{Segalfibrant} for reduced simplicial objects in $\Top_{red}^{\Theta_{n-1}^\op}$ implies that $\delta^*_{\Theta_{n-1}}(f)$ is a trivial Reedy-fibration, cf. Remarks \ref{Segalfibrant2} and \ref{comparison}; it follows then from the commutative triangle (\ref{rs}) that $\sg_{n-1}^*(f)$ is also a realisation weak equivalence and hence, by the inductive hypothesis, a trivial Reedy-fibration. Therefore, since the level-trees $T_k,\,1\leq k\leq m,$ belong to the image of $\sg_{n-1}$, each factor $f(T_k):X(T_k)\to Y(T_k),\,1\leq k\leq m,$ in the lower horizontal map above is a trivial fibration so that the upper horizontal map $f(T):X(T)\to Y(T)$ is a weak equivalence. Consequently, the given realisation-trivial Segal-fibration $f:X\to Y$ is a pointwise weak equivalence and hence a trivial Reedy-fibration as required.

The underlying-space functor $U:\Top^{\Theta_n^\op}\to\Top_*:X\mapsto X(\bar{n})$ is homotopy-right-conservative, since for a Segal-fibrant reduced $\Theta_n$-space $X$, the canonical $n$-graphical decomposition of $T_*$ into a colimit of representable graphs induces a functorial trivial fibration $X(T)\overset{\sim}{\cto}\prod_{\gamma_n(T)}X(\bar{n})$, where there are as many factors as there are vertices in $T$ of height $n$, cf. Section \ref{level-tree} and Remark \ref{factor2}. Finally, for each reduced $\Theta_n$-space $X$, the canonical map $U(X)\to\Omega^n|X|_{\Theta_n}$ factors as $$U(\sg_{n-1}^*X)\overset{\alpha}{\lra}\Omega^{n-1}|\sg_{n-1}^*X|_{\Theta_{n-1}}\overset{\Omega^{n-1}\beta}{\lra}\Omega^{n-1}\Omega|X|_{\Theta_n}.$$If $X$ is cofibrant Segal-fibrant, then so are $\sg_{n-1}^*(X)$ and $\delta^*_{\Theta_{n-1}}(X)$, so that $\alpha$ and $\beta$ are weak equivalences by the inductive hypothesis and by Theorem \ref{Segalfibrant} for reduced simplicial objects in $\Top_{red}^{\Theta_{n-1}^\op}$, cf. Remarks \ref{Segalfibrant2} and \ref{comparison}.\end{proof}

\begin{rmk}\label{comparison}There are two other models for $n$-fold loop spaces (in the sense of Lemma \ref{nfold}) closely related to the preceding one. The first one results from a suitable model structure on \emph{reduced $n$-simplicial spaces}, and can be thought of as a straightforward iteration of the Segal model for $1$-fold loop spaces (an $n$-fold loop space being a loop space in $(n-1)$-fold loop spaces). The second model is due to Bousfield \cite{Bou} and results from a suitable model structure on \emph{$(n-1)$-reduced simplicial spaces}. We are grateful to Rainer Vogt for pointing out to us that the classical diagonal functor $d:\Delta\to\Delta\times\cdots\times\Delta$ cannot be used to relate the reduced $n$-simplicial model to the $(n-1)$-reduced simplicial model since the inverse image functor $d^*$ fails to take reduced $n$-simplicial sets to $(n-1)$-reduced simplicial sets. Instead, we shall use below a realisation functor $\widehat{\Theta}\to\widehat{\Delta}$ which takes reduced $\Theta_n$-sets to $(n-1)$-reduced simplicial sets.

To be precise, an $n$-simplicial space $X(-,\cdots,-)$ is \emph{reduced} if $X([k_1],\dots,[k_n])$ is the one-point space whenever $[k_s]=[0]$ for some $s$. A simplicial space $X(-)$ is \emph{$(n-1)$-reduced} if $X([k])$ is the one-point space for $k<n$. The existence of a suitable localisation of the Reedy model structure on reduced $n$-simplicial spaces follows inductively from Theorem \ref{Segalfibrant} using Remark \ref{Segalfibrant2}. The inductive step is based on the following preservation property of the Segal construction: if $\EE$ is an internal model category having a flat cosimplicial object and fulfilling the realisation and fibrancy axioms (see \ref{Segalfibrant2}), then the localised model category of reduced simplicial $\EE$-objects is again such a category. Bousfield  constructs in \cite{Bou} an explicit suitable localisation of the Reedy model category of $(n-1)$-reduced simplicial spaces.\end{rmk}

\begin{prp}\label{compare1}The diagonal $\delta_n:\Delta^{\times n}\to\Theta_n$ induces a Quillen equivalence between the reduced $\,\Theta_n$- and the reduced $n$-simplicial model for $n$-fold loop spaces.\end{prp}

\begin{proof}By Proposition \ref{rigid}, it suffices to show that $\delta_n$ induces a well-adapted Quillen-adjunction between the models $(|\!-\!|_{\Theta_n},\Top_{red}^{\Theta_n^\op},U)$ and $(|\!-\!|_{\Delta^{\times n}},\Top_{red}^{{\Dop}^{\times n}},U)$ for $n$-fold loop spaces. Observe first that $\delta_n^*$ takes reduced $\Theta_n$-spaces to reduced $n$-simplicial spaces; the compatibility of $\delta^*_n$ with the realisation functors follows then from Proposition \ref{realisation}. The underlying-space of a reduced $n$-simplicial space $X$ is $X([1],\dots,[1])$; the underlying space of the right Kan extension $\delta_{n*}(X)$ is given by the $n$-simplicial ``mapping space'' from $\delta_n^*(\Theta_n[\bar{n}],\partial\Theta_n[\bar{n}])$ to $(X,*)$. The latter is canonically isomorphic to $X([1],\dots,[1])$. Therefore, the adjunction $(\delta_n^*,\delta_{n*})$ is well-adapted. This adjunction is a Quillen adjunction with respect to the Reedy model structures. Since both Reedy structures are localised with respect to realisation weak equivalences, and the left adjoint $\delta_n^*$ is compatible with the realisation functors, the adjunction $(\delta_n^*,\delta_{n*})$ is also a Quillen adjunction after localisation.\end{proof}

\begin{prp}\label{compare2}There is a minimal combinatorial $n$-disk in simplicial sets, whose classifying morphism $d_{n*}:\widehat{\Delta}\rlto\widehat{\Theta}_n:d_n^*$ induces a Quillen equivalence between the reduced $\,\Theta_n$- and the $(n-1)$-reduced simplicial model for $n$-fold loop spaces.\end{prp}

\begin{proof}In the category of simplicial sets, put $D_0=\Delta[0],\,D_1=\Delta[1]$ and for $n>1$, $$D_n=\Delta[n]/(\partial_0\Delta[n-1]\cup\cdots\cup\partial_{n-2}\Delta[n-1]).$$It is then readily verified that\begin{diagram}[noPS,small,silent]D_n&\pile{\lTo^{\partial_n}\\\rTo^{s_{n-1}}\\\lTo^{\partial_{n-1}}}&D_{n-1}&\pile{\lTo^{\partial_{n-1}}\\\rTo^{s_{n-2}}\\\lTo^{\partial_{n-2}}}&D_{n-2}&\cdots&D_2&\pile{\lTo^{\partial_2}\\\rTo^{s_1}\\\lTo^{\partial_1}}&D_1&\pile{\lTo^{\partial_1}\\\rTo^{s_0}\\\lTo^{\partial_0}}&D_0\end{diagram}defines a combinatorial $n$-disk in simplicial sets which realises to the combinatorial $n$-disk underlying the topological $n$-ball, cf. the proofs of Theorem \ref{disk} and Corollary \ref{flat}. Therefore, there is an essentially unique geometric morphism of toposes $d_{n*}:\widehat{\Delta}\rlto\widehat{\Theta}_n:d_n^*$ such that the left adjoint $d_n^*$ takes the generic combinatorial $n$-disk $\Theta[\bar{n}]$ to $D_n$. This left adjoint commutes up to isomorphism with the realisation functors and preserves thus ``boundaries''. In particular, $d_n^*$ takes reduced $\Theta_n$-sets to $(n-1)$-reduced simplicial sets. Since $\Top$ is topological over $\Sets$, the adjunction lifts to an adjunction between reduced $\Theta_n$-spaces and $(n-1)$-reduced simplicial spaces in such a way that the right adjoint commutes up to isomorphism with the underlying-space functors. Moreover, the adjunction is a Quillen adjunction with respect to the canonical Reedy model structures (on the reduced objects) before and after localisation. Proposition \ref{rigid} allows us then to conclude.\end{proof}

\begin{rmk}\label{end}We completely left out the discussion of \emph{operadic} models for $n$-fold loop spaces. There is indeed much to say about these and their relationship to the Segal-type models described above. One of the motivations of this text was to prepare things for a general picture which would incorporate them all. 

An \emph{$E_n$-algebra} (see Boardman-Vogt \cite{BV}, May \cite{M}) models an $n$-fold loop space only if it is \emph{group-complete} so that the classical model structures on $E_n$-algebras have to be localised in order to give rise to models for $n$-fold loop spaces in the sense of Lemma \ref{nfold}. With this caveat in mind, Proposition \ref{rigid} may prove useful in constructing explicit chains of Quillen equivalences relating presheaf-type models and operadic models for $n$-fold loop spaces. 

In this context, let us mention that by \cite[Proposition 1.16]{mi} an \emph{$n$-operad} in Batanin's sense \cite {Ba1} with values in (a cocomplete cartesian closed category) $\EE$ is essentially the same as an $\EE$-enriched category over $\Theta_n$ fulfilling some factorisation properties; this correspondence for $n=1$ has been studied by Thomason \cite{T}. Similarly, a \emph{symmetric operad} with values in $\EE$ is essentially the same as an $\EE$-enriched category over $\Gamma$ fulfilling some factorisation properties; the latter correspondence has been studied by May and Thomason \cite{MT}. It turns out that (under these correspondences) ``base-change'' along the assembly functor $\gamma_n:\Theta_n\to\Gamma$ gives rise to a right Quillen functor from reduced symmetric operads to reduced $n$-operads, the left adjoint of which has the remarkable property that its left derived functor takes the terminal reduced $n$-operad to an $E_n$-operad; this is a reformulation of one of the central results of Batanin \cite{Ba2} and should ultimately lead to a combinatorial characterisation of $E_n$-operads. In particular, combining technics of Thomason \cite{T} and May-Thomason \cite{MT} with Proposition \ref{compare1} should recover some of the comparison results of Dunn \cite{Du} and Fiedorowicz-Vogt \cite{FV}.\end{rmk}

\subsection{The canonical reduced $\Theta_n$-set model for Eilenberg-MacLane spaces}

We come back to the $\Theta_n$-set model $\gamma_n(H\pi)$ of Corollary \ref{EM} for an Eilenberg-MacLane space of type $K(\pi,n)$. For brevity, we shall write $\gamma_n(H\pi)=K(\pi,n)$. The inverse image functors $\delta_n^*:\widehat{\Theta}_n\to\widehat{\Delta^{\times n}}$  and $d_n^*:\widehat{\Theta}_n\to\widehat{\Delta}$ of Propositions \ref{compare1} and \ref{compare2} take $K(\pi,n)$ to the classical $n$-simplicial and simplicial models of $K(\pi,n)$. 

Observe that in \emph{any} model category for $n$-fold loop spaces in the sense of Lemma \ref{nfold}, a \emph{discrete} cofibrant-fibrant object $X$ deloops to an Eilenberg-MacLane space $\Phi(X)$. Indeed, the discrete underlying space $U(X)$ is weakly equivalent to $\Omega^n\Phi(X)$, whence $\Phi(X)$ is $n$-coconnected, and the homotopy-counits of the $\Phi$-$\Psi$-adjunction are $(n-1)$-connected covers, whence $\Phi(X)$ is $(n-1)$-connected; therefore, $\Phi(X)$ is an Eilenberg-MacLane space of type $K(\pi_n(\Phi(X)),n)$; moreover, the isomorphism $\pi_n(\Phi(X))\cong\pi_0(U(X))$ endows the underlying set of $X$ with a canonical group structure. In our case, the delooping functors are given by realisation functors, which explains why the discrete cofibrant-fibrant reduced presheaves are themselves models for Eilenberg-MacLane spaces. We have seen that the inverse image functors $\delta_n^*$ and $d_n^*$ are compatible with the topological realisation functors; nevertheless, the resulting $\CW$-structures depend on the choice of the cell category. We study here the cellular structure of the $\Theta_n$-set $K(\pi,n)$, especially in the case $\pi=\ZZ/2\ZZ$.

By definition, the set $K(\pi,n)(T)$ of $T$-cells of $K(\pi,n)$ is the power-set $\pi^{\gamma_n(T)}$, which can be considered as the set of all labelings of the height-$n$-vertices of $T$ by elements of the group $\pi$. We have to describe the action of the degeneracy operators on these cells. A degeneracy operator $\phi:S\to T$ in $\Theta_n$ corresponds to a subtree-inclusion of $T$ in $S$, cf. Remark \ref{factor2}; the functor $\gamma_n$ takes $\phi:S\to T$ to the projection map $\gamma_n(\phi):\gamma_n(S)\to\gamma_n(T)$ which ``forgets'' about those height-$n$-vertices of $S$ which do not belong to the subtree $T$. Therefore, $K(\pi,n)(\phi):K(\pi,n)(T)\to K(\pi,n)(S)$ labels those height-$n$-vertices of $S$ which do not belong to $T$ by the neutral element of $\pi$. This implies that there are non-degenerate $T$-cells in $K(\pi,n)$ only if $T$ is either of height $0$ (in which case there is a unique non-degenerate $0$-cell) or if $T$ is a \emph{pruned $n$-tree}, i.e. a tree whose input vertices are all of height $n$; in the latter case, the non-degenerate $T$-cells correspond to all labelings of the input vertices by non-neutral elements. 

Recall from Proposition \ref{pregeometric} that the dimension of the cell $|\Theta_n[T]|_{\Theta_n}$ defined by $T$ is the number of edges of $T$, and recall from Remark \ref{factor} that $K(\pi,n)$ realises to a $\CW$-complex with as many cells as there are non-degenerate elements in $K(\pi,n)$. We denote by $f_{n,\pi}^k$ the number of non-degenerate $\pi$-labelled pruned $n$-trees with $n+k$ edges, where a non-degenerate $\pi$-labeling consists of a labeling of the input vertices by non-neutral elements of $\pi$. We just proved

\begin{sprp}\label{cellEM}The realisation of the canonical reduced $\Theta_n$-set $K(\pi,n)$ is a reduced $\CW$-complex with $f_{n,\pi}^k$ cells of dimension $n+k$.\end{sprp}

The number $f_{n,\pi}^k$ only depends on the order of the group $\pi$. In particular, $f_n^k=f_{n,\ZZ/2\ZZ}^k$ counts the number of pruned $n$-trees with $n+k$ edges. This number is a \emph{generalised Fibonacci number} as expressed by the following counting lemma, first noticed by James Dolan, cf. \cite{Do}, and its obvious generalisation below.

\begin{slma}For any $n\geq 1$, if $f_n^k=0$ for $k<0$ and $f_n^0=1$, one has the recursion formula $f_n^k+f_n^{k+1}+\cdots+f_n^{k+n-1}=f_n^{k+n}$. In particular, the sequence $(f_2^k)_{k\geq 0}$ is the classical Fibonacci number sequence.\end{slma}

\begin{proof}Take a pruned $n$-tree with $n+k$ edges and consider the edge-path from its rightmost input vertex to the first vertex with more than one incoming edge. This edge-path has a positive length $m<n$, and the given pruned $n$-tree with $n+k$ edges can be recovered in a unique way from the pruned $n$-tree with $n+k-m$ edges, obtained by removing the above-mentioned edge-path.\end{proof}

\begin{slma}For any $n\geq 1$ and any group $\pi$ of finite order $p$, if $f_{n,\pi}^k=0$ for $k<0$ and $f_{n,\pi}^0=p-1$, one has the recursion formula $(p-1)(f_{n,\pi}^k+f_{n,\pi}^{k+1}+\cdots+f_{n,\pi}^{k+n-1})=f_{n,\pi}^{k+n}$.\end{slma}

The generating function $f_{n,\pi}(t)=\sum_{k\geq 0}f_{n,\pi}^kt^k$ of these generalised, weighted Fibonacci numbers is \emph{rational}. It is readily verified that we have$$f_{n,\pi}(t)=\frac{p-1}{1-(p-1)(t+t^2+\cdots+t^n)}.$$Therefore, the generating function $K(\pi,n)(t)=\sum_{d\geq 0}c_d(|K(\pi,n)|_{\Theta_n})t^d$ for the number $c_d(|K(\pi,n)|_{\Theta_n})$ of $d$-cells of the $\CW$-complex $|K(\pi,n)|_{\Theta_n}$ is rational again:$$K(\pi,n)(t)=1+t^nf_{n,\pi}(t)=\frac{1-(p-1)(t+t^2+\cdots+t^{n-1})}{1-(p-1)(t+t^2+\cdots+t^n)}$$and we get as \emph{virtual Euler-Poincar\'e characteristic} $\chi(K(\pi,n))=K(\pi,n)(-1)$:

\begin{sprp}\label{EP}$\chi(K(\pi,n))=p^{(-1)^n}$ for any group $\pi$ of order $p$.\end{sprp}

This calculation suggests the existence of an Euler-Poincar\'e characteristic for a quite large class of $\CW$-complexes, which would simultaneously be \emph{homotopy-invariant}, \emph{additive} with respect to cellular attachments, and \emph{multiplicative} with respect to fibrations. Proposition \ref{EP} reflects the multiplicativity. It remains however ``mysterious'' why the above-defined characteristic should be an additive homotopy-invariant outside the context of finite $\CW$-complexes. In particular, if under some hypotheses on the class of $\CW$-complexes, the three properties above would hold, then the Euler-Poincar\'e characteristic of a simple $\CW$-complex could be calculated in two completely different ways; either ``cellularly'' like above, or ``cocellularly'', by using the Postnikov decomposition of the space; this second way leads to what Baez \cite{Baez} calls the \emph{homotopy cardinality} of the space. The importance of having an Euler-Poincar\'e characteristic calculable in these two different ways has been emphasised by Loday a long time ago.\vspace{5ex}

\noindent{\sc Universit\'e de Nice, Laboratoire J.-A. Dieudonn\'e, Parc Valrose, 06108 Nice Cedex, France.}\hspace{2em}\emph{E-mail:} cberger$@$math.unice.fr\vspace{2ex}

\end{document}